\documentclass{article}
\usepackage[english]{babel}
\usepackage[letterpaper,top=2cm,bottom=2cm,left=3cm,right=3cm,marginparwidth=1.75cm]{geometry}
\usepackage{amsmath,amsfonts,amsthm,authblk,mathtools,bbm}
\usepackage{graphicx}
\usepackage{tikz}
\usepackage{multirow}
\usepackage{amssymb}

\usepackage[colorlinks=true, allcolors=blue]{hyperref}
\newcommand{\p}{\partial}

\newtheorem{theorem}{Theorem}[section]
\newtheorem{lemma}{Lemma}[section]
\newtheorem{corollary}{Corollary}[section]

\newtheorem{definition}{Definition}[section]
\newtheorem{remark}{Remark}[section]
\newcommand*\diff{\mathop{}\!\mathrm{d}}
\newtheorem*{remark*}{Remark}
\usepackage[pagewise]{lineno}
\usepackage{hyperref}
\usepackage[most]{tcolorbox}

\newtcolorbox{schemebox}{
    colback=gray!3,
    colframe=black,
    boxrule=0.8pt,
    arc=2mm,
    left=2mm,
    right=2mm,
    top=2mm,
    bottom=2mm,
    before skip=10pt,
    after skip=10pt
}

\usepackage{color}

\begin{document}

\title{\normalsize\Large\bfseries Bound-preserving finite volume scheme for a coupled parabolic system with nonlinear and nonlocal drift}


\author[1]{Nicklas J\"averg\r{a}rd\thanks{nicklas.javergard@kau.se}}
\affil[1]{Department of Mathematics and Computer Science, Karlstad University,
Karlstad, Sweden}
\author[1]{Adrian Muntean\thanks{adrian.muntean@kau.se}}

\date{\empty} 

\maketitle

\begin{abstract}
We propose a bound preserving finite-volume scheme for a coupled system of partial differential equations with nonlinear, nonlocal drift. The main contribution of this work is twofold: $(1)$ we show that the proposed scheme preserves physically motivated bounds; $(2)$ we prove that the scheme is well-posed and leads to a sequence of approximate solutions converging to the weak solution of the original problem. This setting is motivated by questions related to morphology formation in organic solar cells. 
\end{abstract}

\section{Introduction}
Let $\Omega\neq \emptyset$ be a parallelepiped in $\mathbb{R}^3$, take $T>0$ to be the final time of the problem, and let $Q_T:=\Omega \times(0, T)$ be the periodic cylinder. We are interested in the following problem: Find $(m, \phi)$ satisfying 

\begin{equation}\label{Eq:MainModel}
\left\lbrace
\begin{aligned}
    &\p_t m = \nabla \cdot \left[\nabla m - 2 \beta (\phi -m^2 ) (\nabla J \star m) \right] \quad\,  \text{in } Q_T,\\
    &\p_t \phi = \nabla \cdot \left[ \nabla \phi - 2 \beta m (1 - \phi) (\nabla J \star m) \right] \qquad\text{in } Q_T,\\
    &m, \phi \quad \Omega\text{-periodic}, \\
    &m(0, x) = m_0(x), \quad \phi(0,x) = \phi_0(x) \qquad\text{for } x\in \bar{\Omega}.
\end{aligned}\right.
\end{equation}
In \eqref{Eq:MainModel}, we use the notation $\star$ to denote the convolution operation, which is defined in \eqref{eq:convolution}.
The system of coupled equations shown in \eqref{Eq:MainModel} describes the evolution of an interacting ternary mixture; the positive values of $m$ refers to one species and the negative values represent another species. The third species is given by $(1-\phi)$. For a more detailed description of this model, see our previous work  \cite{NicklasJ2026} and references cited therein.
\newline
The proposed target system is mathematically interesting as it is strongly coupled through nonlinear and nonlocal drift terms and it can also degenerate at points where $m=0, \phi=1$, or when  $\phi=\pm m$. The presence of the  degeneracies and of the special polynomial structure in the drift nonlinearities offer excellent conditions for producing phase separation within the mixture.  We know that the system \eqref{Eq:MainModel} admits unique (global-in-time) weak solutions (see Theorem 1.1 cf.  \cite{LYONS2024}). We wish to solve \eqref{Eq:MainModel} numerically in a proper way. It is worth mentioning that plenty of numerical analysis and simulation work has been done on the phase separation question for multi-component systems particularly using the Cahn-Hilliard and Allen-Cahn systems of equations concerning bound-preserving schemes. In the multi-component case, it is harder to ensure rigorously that approximation schemes preserve physical bounds.
Preserving physical bounds is important since it facilitates a clearer interpretation of the observed spatial patterns/morphologies. This issue seems to be easier to handle within the context of the finite volume method. For instance, in \cite{RAG2024} an energy-stable bound-preserving finite volume scheme was proposed for a nonlocal Cahn-Hilliard equation. For the Allen-Cahn equation, a finite volume scheme based on operator splitting was presented in \cite{PENG2025}, this scheme was shown to be energy stable and bound preserving. These classical phase-field-type models share many similarities with our system \eqref{Eq:MainModel} in the sense that they all provoke phase separation through degeneracies. The advantage of choosing \eqref{Eq:MainModel} is that is has a rigorous derivation from a stochastic interacting particle system; see \cite{Marra}, for instance. The nonlinear and nonlocal structure of the drifts introduce in our case new mathematical challenges both in constructing a bound-preserving scheme and in proving convergence as we will see throughout this paper.

In this context, our goal is to construct a convergent finite volume scheme approximating the weak solution to \eqref{Eq:MainModel} that preserves the inequalities
\begin{equation}
    0\leq |m(x, t)|\leq \phi(x, t) \leq 1 \qquad \text{for a.e. }(x, t)\in Q_T, \label{eq:2}
\end{equation}
independent of the size of the mesh given that the initial condition satisfies 
\begin{equation*}
0\leq |m_0(x)|\leq \phi_0(x) \leq 1 \qquad \text{for a.e. }x\in \Omega. \label{eq:H}
\end{equation*}
We are interested in using \eqref{Eq:MainModel} to simulate the three dimensional phase separation in ternary mixtures which is a key part of the manufacturing process of organic solar cells. 
Therefore, we want the approximation scheme to preserve the discrete analog of \eqref{eq:2} independently of the size of the mesh whilst still converging to the weak solution to \eqref{Eq:MainModel}. 
The inequalities in \eqref{eq:2} express a volume constraint between the constituents of the mixture, meaning that the three species are not allowed to occupy the same space at the same time.
To simplify the mathematical analysis it is convenient to introduce a new set of unknown functions and balance-equations.
Taking $w:=\phi + m$ and $v:=\phi-m$ leads to the three inequalities $\,v\geq0$, $w\geq 0$, and $w + v \leq 2$ being equivalent to \eqref{eq:2}. Consequently, from now on, we consider the following evolution system: Find $(w, v)$ such that 
\begin{equation}\label{Eq:MainModel_wv}
\left\lbrace
\begin{aligned}
    &\p_t w + \nabla \cdot q_w = 0,
    \qquad q_w:=-\nabla w + 2\beta(1-m)w\nabla J\star m
    \quad\,  \text{in } Q_T,\\
    &\p_t v + \nabla \cdot q_v  = 0,
    \qquad q_v:=-\nabla v -2\beta(1+m)v\nabla J\star m \qquad\,  \text{in } Q_T,\\
    &w,v\quad \Omega\text{-periodic},
    \\
    &w(x, 0)=w_0(x),\qquad
    v(x, 0)=v_0(x) \qquad \text{for }x\in \bar{\Omega}.
\end{aligned}\right.
\end{equation}
In \eqref{Eq:MainModel_wv}, we define $w_0(x) := \phi_0(x) + m_0(x)$ and $v_0(x):=\phi_0(x) - m_0(x)$ for $x\in \bar{\Omega}$. 
For the reader's convenience, we introduce the spaces $H^1_\sharp(\Omega)$ and $H^{-1}_\sharp(\Omega)$. In this context, we denote by \(H^1_\sharp(\Omega)\) the closure of the space of smooth \(\Omega\)-periodic functions with respect to the \(H^1(\Omega)\)-norm. By $H^{-1}_\sharp(\Omega)$ we refer to the dual space of $H^1_\sharp(\Omega)$, i.e. $H^{-1}_\sharp(\Omega) := (H^1_\sharp(\Omega))'$. We refer the reader to \cite{Deimling1985} for details concerning function spaces of periodic functions. \newline
We are interested in approximating the weak solution to the system of equations in \eqref{Eq:MainModel_wv}.
\begin{schemebox}
\begin{definition}\label{Def:MainModelWeak_wv}
    We call a pair $$(w, v) \in \left(L^2(0,T;H^1_\sharp(\Omega))\cap L^\infty(0,T;L^2(\Omega))\right)^2$$ with $$(\p_t w, \p_t v) \in \left( L^2(0,T;H^{-1}_\sharp(\Omega)) \right)^2 $$ a weak solution to \eqref{Eq:MainModel_wv} if the identities,
    \[ \langle  \p_t w , \eta \rangle   + \int_\Omega\nabla w \cdot  \nabla \eta \diff x =  \int_\Omega\beta (1 - m)w(\nabla J \star (w-v)) \cdot \nabla \eta \diff x  , \] and
    \[\langle  \p_t v , \psi \rangle   +\int_\Omega \nabla v \cdot  \nabla \psi \diff x =  -\int_\Omega\beta(m+1)v (\nabla J \star (w-v)) \cdot \nabla \psi \diff x  , \]
    hold for all $(\eta, \psi) \in L^2(0,T;H^1_\sharp(\Omega))\times L^2(0,T;H^1_\sharp(\Omega))$ and $(w(x, 0), v(x, 0))=(w_0(x), v_0(x))$ for $x\in \bar{\Omega}$.
\end{definition}
We refer to finding this pair of weak solutions $(w, v)$ as problem $(P)$.
\end{schemebox}
Our study is organized as follows: In section \ref{sec:Scheme}, we define what we mean by a admissible finite-volume mesh and construct the proposed finite-volume scheme. 
We prove in section \ref{sec:aprioribounds} that the proposed scheme preserves certain physically relevant bounds independent of the size of the mesh and we show energy-type estimates. 
Following this, in section \ref{sec:well-posed}, we show that the fully discrete scheme is well-posed. 
In section \ref{sec:compactness}, we construct a compactness argument via the discrete version of the Lions-Aubin Lemma to show that there exists a convergent subsequence of approximate solutions obtained by solving the scheme in Definition \ref{def:finite-volume-solution} on finer and finer meshes. 
Finally, we show in section \ref{sec:limit} that the limit of such a subsequence satisfies our concept of weak solution presented in Definition \ref{Def:MainModelWeak_wv}.
\section{Construction of the numerical scheme} \label{sec:Scheme}

\subsection{Admissible mesh and assumptions}\label{sec:meshes}
Using the notation from \cite{eymard2000}, we present within this section, what is meant by admissible meshes.
An admissible finite volume mesh of $\Omega$, denoted by $\mathcal{T}$, consists of three families. A family of open disjoint polygonal convex subsets of $\Bar\Omega$ called control volumes and denoted by $K\in\mathcal{T}$. A family of faces denoted $\mathcal{E}$ such that $\mathcal{E}_K$ is the faces associated with control volume $K$. Finally, we introduce a family of points $P=(x_K)_{K\in \mathcal{T}}$. For an admissible mesh the following conditions must be satisfied:
\begin{enumerate}
    \item[$i)$] The closure of the union of all control volumes in $\mathcal{T}$ is $\Bar \Omega$;
    \item[$ii)$] For any $K\in\mathcal{T}$, there exists a subset $\mathcal{E}_K$ of $\mathcal{E}$ such that $\partial K = \Bar K \setminus K = \underset{\sigma \in \mathcal{E}_K}{\bigcup} \Bar \sigma$. Furthermore, $\mathcal{E} = \underset{K\in \mathcal{T}}{\bigcup}\mathcal{E}_K$;

    \item[$iii)$] For any $(K, L) \in \mathcal{T}^2$ with $K\neq L$, either the $2$-dimensional Lebegue measure of $\Bar K \cap \Bar L$ is 0 or $\Bar K \cap \Bar L = \Bar{\sigma}$ for some $\sigma \in \mathcal{E}$, which will be denoted $K|L$;

    \item[$iv)$] The family $P = (x_K)_{K\in \mathcal{T}}$ is such that $x_K\in \Bar K$ (for all $K\in\mathcal{T}$) and, if $\sigma =K|L$, it is assumed that $x_K\neq x_L$, and that the straight line $d_{K,L}$ going through $x_K$ and $x_L$ is orthogonal to $K|L$;

    \item[$v)$] For any $\sigma \in \mathcal{E}$, let $K$ be the control volume such that $\sigma \in \mathcal{E}_K$.
    If $x_K \notin \sigma$, let $d_{K, \sigma}$ be the straight line going through $x_K$ and orthogonal to $\sigma$, assuming that $d_{K, \sigma}\cap \sigma \neq \emptyset$. 
    We denote this point of intersection by $y_\sigma = d_{K, \sigma }\cap \sigma$.
    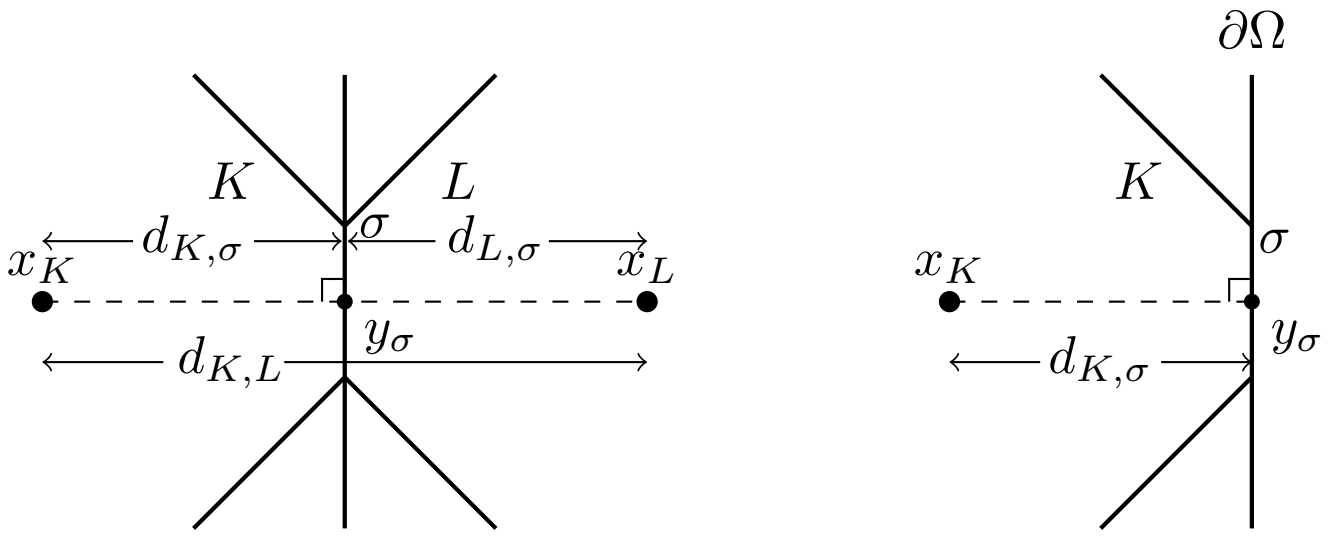
\begin{figure}
        \centering
        \begin{tikzpicture}
            \draw[thick, black] (0,0) -- (1,1);
            \draw[thick, black] (1,1) -- (2,0);
            \draw[thick, black] (1,2) -- (0,3);
            \draw[thick, black] (1,2) -- (2,3);
            \draw[thick, black] (1,3) -- (1,0);
            \fill (-1,1.5) circle (2pt) node[above] {$x_K$};
            \fill (3,1.5) circle (2pt) node[above] {$x_L$};
            \fill (1,1.5) circle (1.5pt) node[below right] {$y_\sigma$};
            \node at (0.25, 1.1) {$d_{K,L}$};
            \node at (0.0, 1.9) {$d_{K,\sigma}$};
            \node at (2.0, 1.9) {$d_{L,\sigma}$};
            \draw[thin, black, <-] (-1, 1.1) --  (-0.2, 1.1);
            \draw[thin, black, ->] (0.6, 1.1) --  (3, 1.1);
            \draw[thin, black, <-](-1, 1.9) --  (-0.4, 1.9);
            \draw[thin, black, ->](0.4, 1.9) --  (0.98, 1.9);
            \draw[thin, black, ->] (2.35, 1.9) --  (3, 1.9);
            \draw[thin, black, <-] (1.02, 1.9) --  (1.6, 1.9);
            \node at (0.25, 2.3) {$K$};
            \node at (1.75, 2.3) {$L$};
            \node at (1.2, 2.0) {$\sigma$};
            \draw[thin, dashed, black] (-1,1.5) -- (3, 1.5);
            \draw[thin, black] (0.85, 1.5) -- (0.85, 1.65) -- (1.0, 1.65);

            \draw[thick, black] (0+6,0) -- (1+6,1);
            \draw[thick, black] (1+6,2) -- (0+6,3);
            \draw[thick, black] (1+6,3) -- (1+6,0);
            \fill (-1+6,1.5) circle (2pt) node[above] {$x_K$};
            \fill (1+6,1.5) circle (1.5pt) node[below right] {$y_\sigma$};
            \draw[thin, dashed, black] (-1+6,1.5) -- (1+6, 1.5);
            
            \node at (6, 1.1) {$d_{K,\sigma}$};
            \draw[thin, black, <-] (5, 1.1) --  (5.6, 1.1);
            \draw[thin, black, ->] (6.4, 1.1) --  (7, 1.1);
            \draw[thin, black] (0.85+6, 1.5) -- (0.85+6, 1.65) -- (1.0+6, 1.65);
            \node at (0.25 + 6, 2.3) {$K$};
            \node at (1.0 + 6, 3.3) {$\partial \Omega$};
            \node at (1.15+6, 1.9) {$\sigma$};
        \end{tikzpicture}
        \caption{Schematic illustration of the notation regarding control volumes, faces, and points for interior (left) and boundary (right) volumes.}
        \label{fig:cells}
    \end{figure}
\end{enumerate}
    We refer the reader to Figure \ref{fig:cells} for an illustration of the entities defined above.\\
    \noindent By $h(\mathcal{T)}:=\sup\{\text{diam} (K):K\in \mathcal{T}\}$ we denote the size of an admissible mesh. We will sometimes refer to $h(\mathcal{T})$ by $h$. To define a numerical scheme for approximating  \eqref{Eq:MainModel_wv} we also need a suitable discretization of the time interval $(0, T)$. We consider a uniform time discretization with the time step $\tau>0$, and set $t^n := n\tau$ for $n\in \{0, 1, \dots , N\}$, where $N\in \mathbb{N}$ satisfying $N\tau =T$.\\
    Furthermore, we assume the following hypotheses to be true:

    \begin{itemize}
        \item[$(H_1)$] $\beta>0$, $J\in C^2_+(\bar{\Omega})$, $J$ is compactly supported and symmetric;
        \item[$(H_2)$] $0\leq |m^0_K|\leq \phi^0_K\leq 1$ for all $K\in\mathcal{T}$;
        \item[$(H_3)$]  There exists $\tau_0 =  \tau_0(h)>0$ such that $\tau \leq \tau_0$ for given $h>0$.
    \end{itemize}

    In $(H_1)$, $\beta >0$ is a natural physical constant. Its meaning is $\frac{1}{T_a}$, where $T_a>0$ is the absolute temperature of the system. For the work within this paper the assumptions on $J$ could be relaxed, however, they are needed in the derivation of the original model.
    The second assumption, $(H_2)$, expresses a volume constraint within each cell $K$ for $n=0$. If $|m^0_K|\to 1$, then this constraint forces $\phi^0_K\to 1$ which means there is no solvent in that cell, recall that the solvent concentration is given by $1-\phi_0(x)$. 
    The last hypothesis, $(H_3)$, is a technical one and it expresses a timestep restriction which is needed to prove uniqueness of solutions to the numerical scheme defined in Definition \ref{def:finite-volume-solution}. \newline
    Before we construct our numerical scheme, we introduce suitable function spaces and norms. To this end, we define
    \begin{equation*}
        \boldsymbol{1}_A(r) :=
        \begin{cases}
            1\quad \text{if }r\in A \\
            0 \quad \text{if }r\notin A
        \end{cases}, \text{ where } \emptyset\neq A\subseteq \Omega.
    \end{equation*}
    Let the space $\mathcal{H_T}$ be the space of piecewise constant functions as  

\begin{equation*}
    \mathcal{H_T} := \bigg\{f:\bar{\Omega}\to \mathbb{R}: \exists(f_K)_{K\in \mathcal{T}}\subset \mathbb{R},\,\,\,\, f(x) =\sum_{K\in\mathcal{T}}f_K\boldsymbol{1}_K(x) \text{ for }x\in\bar{\Omega}\bigg\}.
\end{equation*}
The space $\mathcal{H_T}$ is endowed with the discrete $H^1$-norm given by 

\begin{equation}
    \|f\|^2_{1, \mathcal{T}} := \|f\|^2_{0, \mathcal{T}} + |f|^2_{1, \mathcal{T}},\label{eq:H1norm_def}
\end{equation}
where the two terms on the right-hand side of \eqref{eq:H1norm_def} are defined as:

\begin{align}
    \|f\|^2_{0, \mathcal{T}} &:= \sum_{K\in\mathcal{T}}|K|f_K^2, \\
    |f|^2_{1, \mathcal{T}} &:= \sum_{\sigma\in \mathcal{E}} \frac{|\sigma|}{d_{K,L}}(f_K-f_L)^2.
\end{align}
The discrete $H^{-1}$-norm is the dual norm with respect to the usual $L^2$-scalar product, that is

\begin{equation}
    \|f\|_{-1, \mathcal{T}} := \sup\bigg\{\int_\Omega fg\, dx: g\in \mathcal{H_T}, \|g\|_{1, \mathcal{T}} = 1\bigg\}.\label{eq:H-1norm_def}
\end{equation}
Finally, we define the space $\mathcal{H}_{\mathcal{T},\tau}$ to be the space of piecewise constant in time functions taking values in $\mathcal{H_T}$, namely

\begin{equation} \label{eq:H_T_t}
    \mathcal{H}_{\mathcal{T}, \tau} := \bigg\{f:Q_T\to \mathbb{R} : \exists (f^n)_{n=0}^{N-1},\subset \mathcal{H_T},\,\ f(t, x) = \sum_{n=0}^{N-1}f^n(x)\boldsymbol{1}_{(t^n, t^{n+1})}(t)\bigg\}.
\end{equation}
We equip the space $\mathcal{H}_{\mathcal{T}, \tau}$ with the discrete norm

\begin{equation*}
    \|f\|^2_{1, \mathcal{T}, \tau} := \sum_{n=0}^{N}\tau \|f^{n}\|^2_{1, \mathcal{T}}.
\end{equation*}
Throughout this work, both the model and numerical scheme benefit of periodic boundary conditions. This leads to the following discrete integration by parts formula. For $u\in \mathcal{H_T}$ and a numerical flux $F_{K, \sigma}$ such that $F_{K, \sigma}=-F_{L, \sigma}$ for a face $\sigma=K|L$, we have

\begin{equation}\label{eq:DIBP}
    \sum_{K\in\mathcal{T}}u_K\sum_{\sigma \in \mathcal{E}_k} F_{K, \sigma}= \sum_{\sigma \in \mathcal{E}\\\sigma =K|L} (u_K-u_L)F_{K, \sigma}.
\end{equation}
For an admissible mesh as defined in section \ref{sec:meshes} endowed with periodic boundary conditions, we have 
\begin{equation}
    \sum_{\sigma\in\mathcal{E}_K}|\sigma|d_{K, \sigma} = 3|K|,\label{eq:face-volume-relation}
\end{equation}
where the factor $3$ arising in \eqref{eq:face-volume-relation} stems from posing the problem in three space dimensions.

\subsection{Finite volume discretization of problem $(P)$}\label{sec:Discretize}
Let $(\mathcal{T},\tau)$ be an discretization of the parabolic cylinder $Q_T$ such that $\mathcal{T}$ is an admissible mesh as described in section \ref{sec:meshes}, while $\tau$ represents the uniform timestep. 
We denote by $(w_K^n, v_K^n)$ the discrete unknowns; the values $w^n_K$ and $v^n_K$ are approximations of $w(x_K, nk)$ and $v(x_K, nk)$, respectively. 
We now illustrate the derivation of the scheme for $w$ and simply state the resulting scheme of the same procedure for $v$. 
Integrating \eqref{Eq:MainModel_wv} over an arbitrary control volume $K$ of $\mathcal{T}$ and time interval $(n\tau, (n+1)\tau)$ gives

\begin{equation*}
    \int_{nk}^{(n+1)k}\int_K \partial_t w(x, t) \,dx\, dt - \int_{nk}^{(n+1)k}\int_K \nabla \cdot \big[ \nabla w(x, t) -2\beta (1-m(x, t))w(x, t)\nabla J\star m(x, t)\big]\,dx\, dt =0.
\end{equation*}
Divide this expression by $\tau$ and integrate the first term to obtain

\begin{equation*}
    \frac{1}{\tau} \int_K \big[w(x, t^{n+1}) - w(x, t^n)\big]dx - \frac{1}{\tau}\int_{n\tau}^{(n+1)\tau}\int_K \nabla \cdot \big[ \nabla w(x, t) - 2\beta(1-m(x, t))w(x, t)\nabla J \star m(x, t)\big] \,dx\, dt = 0.
\end{equation*}
Now, we apply Gauss's divergence theorem and choose to approximate the flux by its value at the time step $(n+1)\tau$ in order to arrive at the following implicit time discretization, i.e.

\begin{align}
        \frac{1}{\tau} &\int_K \big[w(x, t^{n+1}) - w(x, t^n)\big]dx \nonumber  \\
        &- \int_{\partial K} \big[ \nabla w(x, t^{n+1}) - 2\beta(1-m(x, t^{n+1}))w(x, t^{n+1})\nabla J \star m(x, t^{n+1})\big] \cdot n_{K} \,d\sigma  = 0,
        \label{eq:WcontinuumApproximation}
\end{align}
where $n_K$ is the outward normal of cell $K$.
As is customary in the finite-volume approach, we approximate functions by their cell average.
Hence, for a function $f\in L^1(K)$, we have 

\begin{equation*}
    f_K := \frac{1}{|K|}\int_K f(x, y, z) \,dx\,dy\,dz.
\end{equation*}
With this notation at hand, we can write:
\begin{equation*}
    \frac{|K|}{\tau}\big(w_K^{n+1} - w_K^n\big) - \int_{\partial K} \big[ \nabla w(x, t^{n+1}) - 2\beta(1-m(x, t^{n+1}))w(x, t^{n+1})\nabla J \star m(x, t^{n+1})\big] \cdot n_{K} \,d\sigma = 0.
\end{equation*}
To move on from \eqref{eq:WcontinuumApproximation}, we need to make choices on how to approximate the flux terms.
To handle the first term we use a standard two-point flux approximation for any cells $K$, $L \in \mathcal{T}$ such that $\sigma =K|L$. This gives 

\begin{equation}\label{eq:twopointfluxapprox}
    -\int _\sigma \nabla w(x, t^{n+1}) \cdot n_K \, d\sigma =  \frac{|\sigma|}{d_{K,L}}(w_K^{n+1} - w_L^{n+1}) +  |\sigma|\mathcal{O}(h),
\end{equation}
which is applicable for all faces of our domain due to the periodic boundary conditions.
For the drift part of the flux, we approximate the face values expressed with subscript $\sigma$ by averages over neighboring cells to get

\begin{equation}\label{eq:fluxapprox}
    \big[2\beta \big(1 - m(x, t^{n+1})\big)w(x, t^{n+1})\nabla J \star m(x, t^{n+1})\big]\cdot n_K \,d\sigma = 2\beta |\sigma|(1-m^{n+1}_\sigma)w_K^{n+1} V^{n+1}_\sigma \cdot n_K + o(h).
\end{equation}
In \eqref{eq:fluxapprox} we have replaced $m(x, t^{n+1})$ by the average of two adjacent cell averages and $w(x, t^{n+1})$ by its cell average as approximations. We expressed the convolution

\begin{equation}
    V(x) := (\nabla J \star m )(x) = \int_\Omega \nabla J (x-y)m(y) \,dy, \label{eq:convolution}
\end{equation}
by its discrete analog
\begin{equation*}
    V^{n+1}_K := \sum_{L\in \mathcal{T}} |L| \nabla J (x_K-x_L)m^{n+1}_L.    
\end{equation*}
Now, we are ready to define the face values $V^{n+1}_\sigma$ and $m^{n+1}_\sigma$ as:

\begin{align*}
    V^{n+1}_\sigma &:= \frac{V^{n+1}_K + V^{n+1}_L}{2},\\
    m^{n+1}_\sigma &:= \frac{m^{n+1}_K+m^{n+1}_L}{2} \quad \text{ for }\sigma = K|L.
\end{align*}
The discrete scheme meant to approximate solutions to problem $(P)$ is as follows:

\begin{schemebox}
\begin{definition}\label{def:finite-volume-solution}
We call $w_{\mathcal{T},\tau}:=(w_K^n)^{n\in\{0,1,\dots, N\}}_{K\in \mathcal{T}}\in \mathcal{H}_{\mathcal{T}, \tau}$ and $v_{\mathcal{T},\tau}:=(v_K^n)^{n\in\{0,1,\dots, N\}}_{K\in \mathcal{T}}\in \mathcal{H}_{\mathcal{T}, \tau}$ solutions to problem $(P_{h, \tau})$ if the following system of equations are satisfied for all $K\in\mathcal{T}$ and for all $n\in\{0,1, \dots, N\}$.
\begin{subequations}\label{Eq:MainModelDisc_wv}
\begin{align}
    &|K| \frac{w_K^{n+1} - w_K^n}{\tau} + \sum_{\sigma\in \mathcal{E}_K} F_{K, \sigma}^{w, n+1} = 0,
    \label{Eq:MainModelDisc_wv_w}
    \\[1ex]
    &|K| \frac{v_K^{n+1} - v_K^n}{\tau} + \sum_{\sigma\in \mathcal{E}_K} F_{K, \sigma}^{v, n+1} = 0,
    \label{Eq:MainModelDisc_wv_v}\\
    &w^0_K := \frac{1}{|K|}\int_K w_0(x), \qquad v^0_K := \frac{1}{|K|}\int_K v_0(x) \quad \text{  for all }K\in \mathcal{T},
\end{align}
\end{subequations}
where 

\begin{subequations}\label{Eq:DiscFlux_wv}

\begin{align}
    F_{K, \sigma}^{w, n+1} :&= \frac{|\sigma|}{d_{K,L}}\big(w_K^{n+1}-w_L^{n+1}\big) + 2\beta|\sigma|\Big(\frac{a^{w,n+1}_{K,\sigma}+\kappa^{n+1}_\sigma}{2}w_K^{n+1} + \frac{a^{w, n+1}_{K,\sigma}-\kappa_\sigma^{n+1}}{2}w_L^{n+1}\Big)\nonumber\\
    &=\frac{|\sigma|}{d_{K,L}}\big(w_K^{n+1}-w_L^{n+1}\big) + 2\beta|\sigma|\Big(a^{w, n+1}_{K, \sigma}\big(\frac{w_K^{n+1}+w_L^{n+1}}{2}\big) + \frac{\kappa_\sigma^{n+1}}{2}\big(w_K^{n+1}-w_L^{n+1}\big)\Big), \label{eq:DiscFlux_w}
    \\[1ex]
    F_{K, \sigma}^{v, n+1} :&= \frac{|\sigma|}{d_{K,L}}\big(v_K^{n+1}-v_L^{n+1}\big) + 2\beta|\sigma|\Big(\frac{a^{v,n+1}_{K,\sigma}+\kappa^{n+1}_\sigma}{2}v_K^{n+1} + \frac{a^{v, n+1}_{K,\sigma}-\kappa_\sigma^{n+1}}{2}v_L^{n+1}\Big)\nonumber\\
    &=\frac{|\sigma|}{d_{K,L}}\big(v_K^{n+1}-v_L^{n+1}\big) + 2\beta|\sigma|\Big(a^{v, n+1}_{K, \sigma}\big(\frac{v_K^{n+1}+v_L^{n+1}}{2}\big) + \frac{\kappa_\sigma^{n+1}}{2}\big(v_K^{n+1}-v_L^{n+1}\big)\Big).\label{eq:DiscFlux_v}
\end{align}
\end{subequations}
In \eqref{eq:DiscFlux_w} and \eqref{eq:DiscFlux_v} we make use of the notation $|\sigma|$ to mean the $2$-dimensional Lebesgue measure of a face between two cells.
The quantities $a^{w, n+1}_{K, \sigma}$, $a^{v, n+1}_{K, \sigma}$ and $\kappa_\sigma^{n+1}$ are given by
\begin{subequations}\label{eq:fluxparams}
\begin{align}
    &a^{w, n+1}_{K, \sigma} := (1-m_\sigma^{n+1}) V^{n+1}_\sigma\cdot n_K,\\
    &a^{v, n+1}_{K, \sigma} := -(1+m_\sigma^{n+1}) V^{n+1}_\sigma\cdot n_K,  \\
    &\kappa^{n+1}_{\sigma} := (1+|m^{n+1}_\sigma|)|V_\sigma^{n+1}\cdot n_K|.
\end{align}
\end{subequations}
\end{definition}
\end{schemebox}

In \eqref{Eq:DiscFlux_wv} we show two equivalent forms of our numerical fluxes, the first line in \eqref{eq:DiscFlux_w} and \eqref{eq:DiscFlux_v} show the upwind structure, since $\kappa_\sigma^{n+1}\geq a^{w, n+1}_{K,\sigma}$ and $\kappa_\sigma^{n+1}\geq a^{v, n+1}_{K,\sigma}$. The second line in \eqref{eq:DiscFlux_w} and \eqref{eq:DiscFlux_v} show that the terms  $\frac{\kappa_\sigma^{n+1}}{2}\big(w_K^{n+1}-w_L^{n+1}\big)$ and $\frac{\kappa_\sigma^{n+1}}{2}\big(v_K^{n+1}-v_L^{n+1}\big)$ have the same structure as the two-point flux approximation in \eqref{eq:twopointfluxapprox}. As we will see later in Lemma \ref{lem:boundedness}, the presence of the additional dissipation is crucial in order to prove the inequality $2\phi_K^{n+1}=w_K^{n+1} + v_K^{n+1}\leq2$.
If we add \eqref{Eq:MainModelDisc_wv_w} and \eqref{Eq:MainModelDisc_wv_v} we arrive at a scheme for $\phi_K^{n}$. The scheme for $\phi^n_K$ looks very similar to the schemes for $w^n_K$ and $v^n_K$.

\begin{equation}\label{eq:MainModelDisc_f}
    \frac{|K|}{\tau}\big(\phi_K^{n+1}-\phi_K^{n}\big) + \sum_{\sigma\in\mathcal{E}_K}F^{\phi, n+1}_{K,\sigma} = 0.
\end{equation}
When we define the numerical flux for $\phi_K^{n+1}$, we must do so in terms of $F^{w, n+1}_{K,\sigma}$ and $F^{v, n+1}_{K,\sigma}$ in order to have a consistent system, as shown in the following calculations

\begin{align}
    \frac{1}{2}\Big(F^{w, n+1}_{K,\sigma} &+ F^{v, n+1}_{K,\sigma}\Big)
    = \frac{|\sigma|}{2d_{K,L}}\big(w_K^{n+1}-w_L^{n+1}\big) + \frac{|\sigma|}{2d_{K,L}}\big(v_K^{n+1}-v_L^{n+1}\big)\nonumber\\
    &+ \beta|\sigma|\Bigg[a^{w,n+1}_{K, \sigma}\Big(\frac{w_K^{n+1}+w_L^{n+1}}{2}\Big) + \frac{\kappa_\sigma^{n+1}}{2}\Big(w_K^{n+1}-w_L^{n+1}\Big)\Bigg] \nonumber \\
    &+\beta|\sigma|\Bigg[a^{v,n+1}_{K, \sigma}\Big(\frac{v_K^{n+1}+v_L^{n+1}}{2}\Big) + \frac{\kappa_\sigma^{n+1}}{2}\Big(v_K^{n+1}-v_L^{n+1}\Big)\Bigg] \nonumber \\
    &=\frac{|\sigma|}{d_{K, L}}\big(\phi_K^{n+1}-\phi_L^{n+1}\big) + \frac{2\beta|\sigma|\kappa_\sigma^{n+1}}{2}\big(\phi_K^{n+1}-\phi_L^{n+1}\big) + \nonumber \\
    &+\beta|\sigma|\Bigg[(1-m_\sigma^{n+1})V_\sigma^{n+1}\cdot n_K\big(\frac{w_K^{n+1}+w_L^{n+1}}{2}\big) - (1+m_\sigma^{n+1})V_\sigma^{n+1}\cdot n_K\big(\frac{v_K^{n+1}+v_L^{n+1}}{2}\big)\Bigg]. \nonumber 
\end{align}
Collecting all terms with $m_\sigma^{n+1}$ and utilizing $2m_K = w_K-v_K$ and $2\phi_K = w_K + v_K$ for all $K\in \mathcal{T}$ we get
\begin{align}
     \frac{1}{2}\Big(F^{w, n+1}_{K,\sigma} + F^{v, n+1}_{K,\sigma}\Big)&=\frac{|\sigma|}{d_{K, L}}\big(\phi_K^{n+1}-\phi_L^{n+1}\big) + \frac{2\beta|\sigma|\kappa_\sigma^{n+1}}{2}\big(\phi_K^{n+1}-\phi_L^{n+1}\big) + \nonumber \\
        &+\beta|\sigma|\bigg[\big(\frac{w_K^{n+1}-v_K^{n+1}}{2} + \frac{w_L^{n+1}-v_L^{n+1}}{2}\big)V_\sigma^{n+1}\cdot n_K\bigg] \nonumber\\ 
        &-\beta |\sigma|\bigg[m_\sigma^{n+1}\big(\frac{w_K^{n+1}+v_K^{n+1}}{2} + \frac{w_L^{n+1}+v_L^{n+1}}{2}\big)V_\sigma^{n+1}\cdot n_K\bigg] \nonumber \\
    &=\frac{|\sigma|}{d_{K, L}}\big(\phi_K^{n+1}-\phi_L^{n+1}\big) + \frac{2\beta|\sigma|\kappa_\sigma^{n+1}}{2}\big(\phi_K^{n+1}-\phi_L^{n+1}\big) + \nonumber \\
    &+\beta|\sigma|\Bigg[\big(m_K^{n+1} + m_L^{n+1}\big)V_\sigma^{n+1}\cdot n_K-m_\sigma^{n+1}\big(\phi_K^{n+1} + \phi_L^{n+1}\big)V_\sigma^{n+1}\cdot n_K\Bigg] \nonumber \\
    &=\frac{|\sigma|}{d_{K, L}}\big(\phi_K^{n+1}-\phi_L^{n+1}\big) + \frac{2\beta|\sigma|\kappa_\sigma^{n+1}}{2}\big(\phi_K^{n+1}-\phi_L^{n+1}\big) \nonumber \\
    &\quad+ \beta|\sigma|\Bigg[m_\sigma^{n+1}\big(2-\phi_K^{n+1} - \phi_L^{n+1}\big)V_\sigma^{n+1}\cdot n_K\Bigg] \nonumber \\
    &=\frac{|\sigma|}{d_{K, L}}\big(\phi_K^{n+1}-\phi_L^{n+1}\big) \nonumber \\
    &\quad+ 2\beta|\sigma|\Bigg[\frac{\kappa_\sigma^{n+1}}{2}\big(\phi_K^{n+1}-\phi_L^{n+1}\big) + \frac{m_\sigma^{n+1}}{2}\big(2 - \phi_K^{n+1} - \phi_L^{n+1}\big)V_\sigma^{n+1}\cdot n_K\Bigg] \nonumber \\
    &=\frac{|\sigma|}{d_{K, L}}\big(\phi_K^{n+1}-\phi_L^{n+1}\big) \nonumber \\
    &\quad + 2\beta|\sigma|\Bigg[\frac{\kappa_\sigma^{n+1}}{2}\big(-(1-\phi_K^{n+1})+(1-\phi_L^{n+1})\big) \nonumber \\
    &\quad+ \frac{m_\sigma^{n+1}}{2}\big(1 - \phi_K^{n+1}+ 1 - \phi_L^{n+1}\big)V_\sigma^{n+1}\cdot n_K\Bigg].\nonumber 
\end{align}
Finally, we define $a^{\phi, n+1}_{K,\sigma}:=m_\sigma^{n+1}V_\sigma^{n+1}\cdot n_K$ and  

\begin{equation}\label{eq:DiscFlux_f}
    F^{\phi, n+1}_{K, \sigma}:=\frac{|\sigma|}{d_{K,L}}\big(\phi_K^{n+1}-\phi_L^{n+1} \big) + 2\beta |\sigma|\bigg[\frac{a^{\phi, n+1}_{K, \sigma}-\kappa_\sigma^{n+1}}{2}\big(1-\phi_K^{n+1}\big) + \frac{a^{\phi, n+1}_{K, \sigma}+\kappa_\sigma^{n+1}}{2}\big(1-\phi_L^{n+1}\big)\bigg].
\end{equation}
Now, we have all ingredients which we need to prove the desired bounds stated in the next section.

\section{\emph{A priori} bounds of solutions to $(P_{h, \tau})$}\label{sec:aprioribounds}
In this section we prove a number of \emph{a priori} bounds that facilitates further the mathematical analysis of the approximate solutions to $(P_{h, \tau})$ and the convergence of said solutions to the weak form of \eqref{Eq:MainModel_wv}.
\noindent Within this section, we use $\kappa_\sigma^{n+1}\ge \max\Big\{|a^{w, n+1}_{K,\sigma}|, |a^{v, n+1}_{K, \sigma}|\Big\}$. We also use the fact that our scheme is locally conservative, i.e. $F^{w, n+1}_{K, \sigma}=-F^{w, n+1}_{L, \sigma}$ and $F^{v, n+1}_{K, \sigma}=-F^{v, n+1}_{L, \sigma}$. These ingredients are needed to  prove the positivity and boundedness of solutions to problem $(P_{h, \tau})$.

\begin{lemma}[Positivity of $w_K^n$ and $v_K^n$, $n\in\mathbb
N$]\label{lem:positivity}

    Assume $(H_1)$ and $(H_2)$. Then for every $n\in\mathbb{N}$, the solution $(w_K^n, v_K^n)$ to problem $(P_{h, \tau})$ in the sense of Definition \ref{def:finite-volume-solution} is positive for all $K\in \mathcal{T}$.
\end{lemma}
\begin{proof}
    The proof proceeds by induction. By our assumptions, $(w^0_K, v_K^0)$ is component-wise positive for all $K\in \mathcal{T}$. As induction hypothesis, we assume that $(w^n_K, v_K^n)$ is positive. We now wish to prove that $(w^{n+1}_K, v_K^{n+1})$ is positive component-wise. \newline
    We prove the positivity explicitly for $w^{n+1}_K$ and point out that exactly the same argument also applies to $v^{n+1}_K$.
    We introduce the set $$\mathcal{N} := \{K\in\mathcal{T} : w_K^{n+1} < 0\},$$
    and assume that $\mathcal{N}\neq \emptyset$. We sum \eqref{Eq:MainModelDisc_wv_w} over all $K\in \mathcal{N}$ to get 
    
    \begin{equation}\label{eq:subsum_w}
        \sum_{K\in \mathcal{N}}|K|\frac{w^{n+1}_K-w^n_K}{\tau} + \sum_{K\in \mathcal{N}}\sum_{\sigma\in \mathcal{E}_K}F^{w, n+1}_{K, \sigma}=0.
    \end{equation}
    For any $K\in \mathcal{N}$, we have $w^{n+1}_K<0$ and $w^n_K\geq 0$ which leads to
\begin{equation}\label{eq:lem31_1}
\frac{|K|}{\tau}(w^{n+1}_K-w_K^n)< 0.
\end{equation}
We now turn our attention to the structure of the flux. We deal with the flux in two cases:
\begin{itemize}
    \item[$i)$] $w^{n+1}_K$ for $K\in \mathcal{N}$ and $w^{n+1}_L$ for $L \in \mathcal{N}$;
    \item[$ii)$] $w^{n+1}_K$ for $K\in \mathcal{N}$ and $w^{n+1}_L$ for $L \notin \mathcal{N}$.
\end{itemize}
Concerning case $i)$, it is clear that the fluxes will cancel when both $K$ and $L$ lie in $\mathcal{N}$ since $F^{w, n+1}_{K,\sigma} = - F^{w, n+1}_{L, \sigma}$ for $\sigma = K|L.$ In the case $ii)$, we have $w^{n+1}_K<0$ and $w^{n+1}_L\geq 0$, which leads to 
\begin{equation}\label{eq:lem31_2}
    \frac{|\sigma|}{d_{K,L}}(w^{n+1}_K - w_L^{n+1})<0,
\end{equation}
and
\begin{equation}\label{eq:lem31_3}
    \beta|\sigma|\kappa_\sigma^{n+1}(w^{n+1}_K - w_L^{n+1})<0. 
\end{equation}
What concerns the last term, we know $a^{w,n+1}_{K, \sigma} + \kappa_\sigma^{n+1}\geq 0$ and $a^{w,n}_{K, \sigma}-\kappa_\sigma^{n+1}\leq 0$. Thus we can infer 
\begin{align}
    \frac{a^{w,n+1}_{K, \sigma}+\kappa_\sigma^{n+1}}{2}w^{n+1}_K &\leq 0, \label{eq:lem31_4}\\
    \frac{a^{w,n+1}_{K, \sigma}-\kappa_\sigma^{n+1}}{2}w^{n+1}_L &\leq 0. \label{eq:lem31_5}
\end{align}
Collecting what we know from \eqref{eq:lem31_1}, \eqref{eq:lem31_2}, \eqref{eq:lem31_3}, \eqref{eq:lem31_4} and \eqref{eq:lem31_5}, we get

\begin{align} 
    \sum_{K\in \mathcal{N}}|K|&\frac{w^{n+1}_K-w^n_K}{\tau} + \sum_{K\in \mathcal{N}}\sum_{\sigma\in \mathcal{E}_K}\frac{|\sigma|}{d_{K,L}}(w^{n+1}_K-w^{n+1}_L) + \nonumber \\
    &\sum_{K\in \mathcal{N}}\sum_{\sigma\in\mathcal{E}_K}2\beta |\sigma|\bigg(a^{w, n+1}_{K, \sigma}\big(\frac{w_K^{n+1} + w_L^{n+1}}{2}\big) + \frac{\kappa_\sigma^{n+1}}{2}\big(w_K^{n+1}-w_L^{n+1}\big)\bigg) < 0.\label{eq:lem_contradiction}
\end{align}
Comparing \eqref{eq:lem_contradiction} with \eqref{eq:subsum_w}, we arrive at a contradiction. Therefore, we conclude that $\mathcal{N}=\emptyset$, which proves the result for $w_K^{n+1}$. The proof of positivity of $v_K^{n+1}$ follows the same arguments. This can be simply seen by observing that the upwind structure of the equations building up the approximating sequences $(w_K^{n})$ and $(v_K^{n})$ are is identical\footnote{Look at the sign of each term in the proposed scheme.}. Therefore, we omit to show the full details for $v_K^{n+1}$.
\end{proof}

\begin{lemma}[Uniform upper bound for $\phi_K^{n}, n \in \mathbb{N}$]\label{lem:boundedness}
     Assume $(H_1)$ and $(H_2)$. Then for every $n\in\mathbb{N}$, the solution $(w_K^n, v_K^n)$ to problem $(P_{h, \tau})$ in the sense of definition \ref{def:finite-volume-solution} satisfies the bound $$\frac{w_K^n + v_K^n}{2} = \phi_K^n\leq 1$$ for all $K\in \mathcal{T}$.
\end{lemma}
\begin{proof}
     The proof proceeds by induction. As a starting point, we know that $\phi^0_K\leq 1$ component-wise. We assume by induction that $\phi^n_K\leq 1$ for an arbitrary $n$. We wish now to prove that this implies $\phi^{n+1}_K\leq 1$. \newline
    We introduce the set

    \begin{equation*}
        \mathcal{P} := \{K\in\mathcal{T}:\phi_K^{n+1}>1\},
    \end{equation*}
    and assume $\mathcal{P}\neq\emptyset$. We sum \eqref{eq:MainModelDisc_f} over all $K\in\mathcal{P}$ to get
    
    \begin{equation}
        \sum_{K\in\mathcal{P}}|K|\frac{\phi_K^{n+1}-\phi_K^{n}}{\tau} + \sum_{K\in\mathcal{P}}\sum_{\sigma\in\mathcal{E}_K}F^{\phi, n+1}_{K, \sigma} = 0.\label{eq:subsum_f}
    \end{equation}
    For any $K\in\mathcal{P}$, we have $\phi_K^{n+1}>1$ and $\phi_K^n\leq 1$. This yields
    
    \begin{equation}
        \frac{|K|}{\tau}\big(\phi_K^{n+1}-\phi_K^{n}\big)>0.\label{eq:lem32_1}
    \end{equation}
    Turning our attention to the structure of the flux, we encounter two separate cases:

    \begin{itemize}
        \item[$i)$] $\phi^{n+1}_K$ for $K\in \mathcal{P}$ and $\phi^{n+1}_L$ for $L \in \mathcal{P}$;
        \item[$ii)$] $\phi^{n+1}_K$ for $K\in \mathcal{P}$ and $\phi^{n+1}_L$ for $L \notin \mathcal{P}$.
    \end{itemize}
    When addressing case $i)$, we see that the fluxes will cancel if both $K$ and $L$ are in $\mathcal{P}$ since $F^{\phi, n+1}_{K, \sigma}= -F^{\phi, n+1}_{L, \sigma}$ for $\sigma=K|L$. What concerns the second case $ii)$, we have $\phi^{n+1}_K>1$ and $\phi^{n+1}_L\leq 1$. This leads to 

    \begin{equation}
        \frac{|\sigma|}{d_{K,L}}\big(\phi^{n+1}_K-\phi^{n+1}_L\big)>0, \label{eq:lem32_2}
    \end{equation}
    and 

    \begin{equation}
         \beta|\sigma|\kappa_\sigma^{n+1}\big(\phi^{n+1}_K-\phi^{n+1}_L\big)>0.\label{eq:lem32_3}
    \end{equation}
    Regarding the last two terms, we know $a^{\phi,n+1}_{K, \sigma} - \kappa_\sigma^{n+1}\leq 0$ and $a^{\phi,n+1}_{K, \sigma} + \kappa_\sigma^{n+1}\geq 0$. Therefore

    \begin{align}
        \frac{a^{\phi,n+1}_{K, \sigma}-\kappa_\sigma^{n+1}}{2}\big(1-\phi^{n+1}_K\big) &\geq 0, \label{eq:lem32_4}\\
        \frac{a^{\phi,n}_{K, \sigma}+\kappa_\sigma^{n+1}}{2}\big(1-\phi^{n+1}_L) &\geq 0. \label{eq:lem32_5}
    \end{align}
    Collecting the results in \eqref{eq:lem32_1}, \eqref{eq:lem32_2}, \eqref{eq:lem32_3}, \eqref{eq:lem32_4}, and \eqref{eq:lem32_5}, we get 

    \begin{align} 
        \sum_{K\in \mathcal{P}}|K|&\frac{\phi^{n+1}_K-\phi^n_K}{\tau} + \sum_{K\in \mathcal{P}}\sum_{\sigma\in \mathcal{E}_K}\frac{|\sigma|}{d_{K,L}}(\phi^{n+1}_K-\phi^{n+1}_L) + \nonumber \\
        &\sum_{K\in \mathcal{P}}\sum_{\sigma\in\mathcal{E}_K}2\beta |\sigma|\bigg(\frac{a^{\phi, n+1}_{K, \sigma}-\kappa_\sigma^{n01}}{2}\big(1 - \phi_K^{n+1}\big) + \frac{a^{\phi, n+1}_{K,\sigma}-\kappa_\sigma^{n+1}}{2}\big(1-\phi_L^{n+1}\big)\bigg) > 0.\label{eq:lem32_contradiction}
    \end{align}
    Comparing \eqref{eq:lem32_contradiction} with \eqref{eq:subsum_f}, we arrive at a contradiction. Therefore, we conclude that $\mathcal{P}=\emptyset$, which finally proves that $\phi_K^{n+1}\leq 1$. 
\end{proof}

We point out that the upper bound in Lemma \ref{lem:boundedness} could be obtained this easily because of the artificial diffusion added in \eqref{eq:DiscFlux_w} and \eqref{eq:DiscFlux_v}, which allowed us to write \eqref{eq:DiscFlux_f} with the proper upwind structure.
\begin{remark}
    Combining the conclusion of Lemma \ref{lem:positivity} and Lemma \ref{lem:boundedness}, we obtain that 
    \begin{equation*}
        w^n_K\geq 0, \qquad v^n_K\geq 0, \quad w^n_K+v^n_K\leq 2,
    \end{equation*}
    or equivalently, the inequality
    \begin{equation*}
        0\leq |m^n_K|\leq \phi^n_K\leq 1
    \end{equation*}
    holds true for all $n\in \{0, 1, \dots, N\}$ and $K\in \mathcal{T}$.
\end{remark}
Next we give another upper bound for $w^n_K$ and $v^n_K$ that stems directly from Lemma \ref{lem:positivity} and the conservative nature of the flux. 
\begin{lemma}[Boundedness of $(w_K^n, v_K^n), n\in \mathbb{N}$]\label{lem:aprioribound}
    Assume $(H_1)$ and $(H_2)$. Then there exists a constant $M>0$ such that $w^0_K\leq M$ and $v^0_K\leq M$ for all $K\in\mathcal{T}$. Then for every fixed $h:=\underset{K\in\mathcal{T}}{\max}\{|K|\}$ there exists a constant $M_h>0$ (depending on $M$ and $h$) such that $$w_K^n\leq M_h, \qquad v_K^n\leq M_h$$ for all $n\in\{0, 1, \dots, N\}$ and all $K\in\mathcal{T}$.
\end{lemma}
\begin{proof}
    Summing up \eqref{Eq:MainModelDisc_wv} over each $K\in \mathcal{T}$, we get
    \begin{align*}
         &\sum_{K\in \mathcal{T}}|K|\frac{w^{n+1}_K-w^n_K}{\tau} + \sum_{K\in\mathcal{T}}\sum_{\sigma\in \mathcal{E}_K}F^{w, n+1}_{K, \sigma}=0, \\
         &\sum_{K\in \mathcal{T}}|K|\frac{v^{n+1}_K-v^n_K}{\tau} + \sum_{K\in\mathcal{T}}\sum_{\sigma\in \mathcal{E}_K}F^{v, n+1}_{K, \sigma}=0. 
    \end{align*}
    Due to the periodic boundary conditions and the fact that the flux is conservative, we have \newline $\underset{K\in\mathcal{T}}{\sum}\underset{\sigma\in\mathcal{E}_K}{\sum}F^{\cdot, n+1}_{K, \sigma}=0$ and consequently, we know
    \begin{align*}
        &\sum_{K\in \mathcal{T}}|K|\frac{w^{n+1}_K-w^n_K}{\tau} = 0, \\
        &\sum_{K\in \mathcal{T}}|K|\frac{w^{n+1}_K}{\tau} = \sum_{L\in \mathcal{T}}|L|\frac{w^n_L}{\tau}.
    \end{align*}
    Every term in the last two sums is non-negative by Lemma \ref{lem:positivity}. Therefore, we may write 
    \begin{equation*}
        |K|w^{n+1}_K \leq \sum_{L\in\mathcal{T}}|L|w^{n+1}_L = \sum_{L\in\mathcal{T}}|L|w^0_L\leq M\sum_{L\in\mathcal{T}}|L|=M|\Omega|.
    \end{equation*}
    Taking
    \begin{equation}\label{eq:Mh}
        M_h := \frac{M|\Omega|}{|K|}
    \end{equation}
    we obtain the desired result for $w_K^{n}$. The proof for $v_K^n$ follows the same steps and is therefore omitted.
\end{proof}
\begin{remark}
    In proving Lemma \ref{lem:aprioribound} we have only used the positivity of solutions shown in Lemma \ref{lem:positivity} and the conservation of mass throughout the system. Lemma \ref{lem:aprioribound} will be used in section \ref{sec:existence} when we apply Brouwer's fixed-point theorem to show that at least one solution exists to problem $(P_{h,\tau})$.
\end{remark}

\begin{lemma}[Energy-type estimates]\label{lem:EnergyBound}
    Assume $(H_1)$-$(H_3)$. Let $w_{\mathcal T,\tau},v_{\mathcal T,\tau}\in\mathcal{H}_{\mathcal{T}, \tau}$ be approximate solutions generated via problem $(P_{h, \tau})$. Then there exists a constant $C>0$, which is independent of $h$ and $\tau$, such that
    \begin{align}
        \sup_{0\leq n \leq N} \bigg\{\|w^n\|^2_{0, \mathcal{T}} + \|v^n\|^2_{0, \mathcal{T}}\bigg\} \leq C, \label{eq:BochnerL2_functions} \\
        \sum_{n=0}^{N}\tau \bigg(|w^{n+1}|^2_{1, \mathcal{T}} + |v^{n+1}|^2_{1, \mathcal{T}} \bigg) \leq C, \label{eq:BochnerL2_gradient}
    \end{align}
    and 
    \begin{equation}
        \sum_{n=0}^{N-1}\tau \bigg(\bigg\|\frac{w^{n+1} - w^{n}}{\tau}\bigg\|^2_{-1, \mathcal{T}} + \bigg\|\frac{v^{n+1} - v^{n}}{\tau}\bigg\|^2_{-1, \mathcal{T}}\bigg) \leq C. \label{eq:Bochner_H-1_dt}
    \end{equation}
    Consequently, it also holds

    \begin{equation}
        \|w_{\mathcal{T}, \tau}\|_{1, \mathcal{T}, \tau} + \|v_{\mathcal{T}, \tau}\|_{1, \mathcal{T}, \tau} \leq C. \label{eq:energybound}
    \end{equation}
\end{lemma}

\begin{proof}
We restate the equations that we will use to obtain the necessary bounds in the interest for legibility.

\begin{subequations}
\begin{align*}
    &|K| \frac{w_K^{n+1} - w_K^n}{\tau} + \sum_{\sigma\in \mathcal{E}_K} \bigg\{\frac{\sigma}{d_{K,L}}(w^{n+1}_K-w^{n+1}_L) + \beta|\sigma|\big(\kappa^{n+1}_\sigma(w^{n+1}_K-w^{n+1}_L)  +a^{w, n+1}_{K, \sigma}(w^{n+1}_K + w^{n+1}_L) \big)\bigg\} = 0,
    \\[1ex]
    &|K| \frac{v_K^{n+1} - v_K^n}{\tau} + \sum_{\sigma\in \mathcal{E}_K} \bigg\{\frac{\sigma}{d_{K,L}}(v^{n+1}_K-v^{n+1}_L) + \beta|\sigma|\big(\kappa^{n+1}_\sigma(v^{n+1}_K-v^{n+1}_L ) +a^{v, n+1}_{K, \sigma}(v^{n+1}_K + v^{n+1}_L)\big)\bigg\} = 0.
\end{align*}
\end{subequations}
To obtain the wanted estimates we multiply the first and the second equation by $w^{n+1}_K$ and $v_K^{n+1}$, respectively. We start by proving   \eqref{eq:BochnerL2_functions} and \eqref{eq:BochnerL2_gradient}. We will deal with \eqref{eq:Bochner_H-1_dt} separately. 
\paragraph{Proof of \eqref{eq:BochnerL2_functions} and \eqref{eq:BochnerL2_gradient}}
We begin by rewriting the time derivative in a slightly different form. In the following expression

\begin{align*}
    &\sum_{K\in \mathcal{T}}\frac{|K|}{\tau} (w^{n+1}_K-w^{n}_K)w^{n+1}_K, \\
    &\sum_{K\in \mathcal{T}}\frac{|K|}{\tau} (v^{n+1}_K-v^{n}_K)v^{n+1}_K,
\end{align*}
we will use 
\begin{equation*}
    (a-b)a = a^2 -ab = a^2 -ab + \frac{b^2}{2} - \frac{b^2}{2} = \frac{1}{2}(a-b)^2 + \frac{1}{2}(a^2-b^2)
\end{equation*}
to obtain

\begin{align}
    \sum_{K\in \mathcal{T}}\frac{|K|}{\tau} (w^{n+1}_K-w^{n}_K)w^{n+1}_K &= \frac{1}{2\tau}\Bigg(\sum_{K \in \mathcal{T}}\bigg\{|K|(w^{n+1}_K - w^{n}_K)^2 + |K|\big((w^{n+1}_K)^2 - (w^n_K)^2\big)\bigg\}\Bigg), \label{eq:dwdt_estimate}\\ 
    \sum_{K\in \mathcal{T}}\frac{|K|}{\tau} (v^{n+1}_K-v^{n}_K)v^{n+1}_K &= \frac{1}{2\tau}\Bigg(\sum_{K \in \mathcal{T}}\bigg\{|K|(v^{n+1}_K - v^{n}_K)^2 + |K|\big((v^{n+1}_K)^2 - (v^n_K)^2\big)\bigg\}\Bigg).\label{eq:dvdt_estimate}
\end{align}
Turning our attention to the diffusion terms, we rewrite them as follows

\begin{align}
    \sum_{K\in \mathcal{T}}\sum_{\sigma\in\mathcal{E}_K}\frac{|\sigma|}{d_{K,L}}(w^{n+1}_K - w^{n+1}_{L})w^{n+1}_K &= \sum_{\sigma\in\mathcal{E}} \frac{|\sigma|}{d_{K,L}}\bigg\{(w^{n+1}_K - w^{n+1}_{L})w^{n+1}_K + (w^{n+1}_L - w^{n+1}_{K})w^{n+1}_L\bigg\} \nonumber \\
    &=\sum_{\sigma\in\mathcal{E}} \frac{|\sigma|}{d_{K,L}}(w^{n+1}_K - w^{n+1}_{L})^2 = |w^{n+1}|^2_{1, \mathcal{T}}, \label{eq:grad_w_estimate} \\
    \sum_{K\in \mathcal{T}}\sum_{\sigma\in\mathcal{E}_K}\frac{|\sigma|}{d_{K,L}}(v^{n+1}_K - v^{n+1}_{L})v^{n+1}_K &= \sum_{\sigma\in\mathcal{E}} \frac{|\sigma|}{d_{K,L}}\bigg\{(v^{n+1}_K - v^{n+1}_{L})v^{n+1}_K + (v^{n+1}_L - v^{n+1}_{K})v^{n+1}_L\bigg\} \nonumber \\
    &=\sum_{\sigma\in\mathcal{E}} \frac{|\sigma|}{d_{K,L}}(v^{n+1}_K - v^{n+1}_{L})^2 = |v^{n+1}|^2_{1, \mathcal{T}}. \label{eq:grad_v_estimate} 
\end{align}
To obtain, \eqref{eq:grad_w_estimate}-\eqref{eq:grad_v_estimate}, we have applied the discrete integration-by-parts formula in \eqref{eq:DIBP}. 
We will apply the same idea to deal with the artificial diffusion terms. The result will be non-negative and can later be disregarded without altering the resulting inequality.

\begin{align}
    \sum_{K\in \mathcal{T}}\sum_{\sigma\in\mathcal{E}_K}\beta |\sigma| \kappa_\sigma^{n+1} (w^{n+1}_K - w^{n+1}_{L})w^{n+1}_K &= \sum_{\sigma\in\mathcal{E}} \beta |\sigma|\kappa_\sigma^{n+1} \bigg\{(w^{n+1}_K - w^{n+1}_{L})w^{n+1}_K + (w^{n+1}_L - w^{n+1}_{K})w^{n+1}_L\bigg\} \nonumber \\
    &=\sum_{\sigma\in\mathcal{E}} \beta |\sigma|\kappa_\sigma^{n+1} (w^{n+1}_K - w^{n+1}_{L})^2,  \label{eq:kappa_w_estimate} \\
    \sum_{K\in \mathcal{T}}\sum_{\sigma\in\mathcal{E}_K}\beta |\sigma|\kappa_\sigma^{n+1}(v^{n+1}_K - v^{n+1}_{L})v^{n+1}_K &= \sum_{\sigma\in\mathcal{E}} \beta |\sigma|\kappa_\sigma^{n+1}\bigg\{(v^{n+1}_K - v^{n+1}_{L})v^{n+1}_K + (v^{n+1}_L - v^{n+1}_{K})v^{n+1}_L\bigg\} \nonumber \\
    &=\sum_{\sigma\in\mathcal{E}} \beta |\sigma|\kappa_\sigma^{n+1}(v^{n+1}_K - v^{n+1}_{L})^2. \label{eq:kappa_v_estimate} 
\end{align}
To estimate the drift term, we rewrite the sum in the same way as when we estimated the diffusive terms. We also make use of $a^{w, n+1}_{K, \sigma} = -a^{w, n+1}_{L, \sigma}$ and $a^{v, n+1}_{K, \sigma} = - a^{v, n+1}_{L, \sigma}$. This yields 

\begin{align}
      \sum_{K\in \mathcal{T}}\sum_{\sigma\in\mathcal{E}_K}\beta |\sigma|a^{w, n+1}_{K, \sigma}(w^{n+1}_{K} &+ w^{n+1}_{L})w^{n+1}_{K} = \sum_{\sigma\in\mathcal{E}}\beta |\sigma|\bigg\{a^{w, n+1}_{K, \sigma}(w^{n+1}_{K} + w^{n+1}_{L})w^{n+1}_{K} + a^{w, n+1}_{L, \sigma}(w^{n+1}_{K} + w^{n+1}_{L})w^{n+1}_{L} \bigg\} \nonumber \\
      &=\sum_{\sigma\in\mathcal{E}}\beta |\sigma|a^{w, n+1}_{K, \sigma}\bigg\{(w^{n+1}_{K} + w^{n+1}_{L})w^{n+1}_{K} -(w^{n+1}_{K} + w^{n+1}_{L})w^{n+1}_{L} \bigg\} \nonumber \\
      &=\sum_{\sigma\in\mathcal{E}}\beta |\sigma|a^{w, n+1}_{K, \sigma}\bigg(\big(w^{n+1}_K\big)^2 - \big(w^{n+1}_L\big)^2\bigg) \nonumber \\
      &\leq \sum_{\sigma\in\mathcal{E}}\beta |\sigma|\, |a^{w, n+1}_{K, \sigma}||w^{n+1}_{K} + w^{n+1}_{L}||w^{n+1}_{K}-w^{n+1}_{L}|\sqrt\frac{d_{K, L}}{d_{K,L}}  \nonumber \\
      &\leq \sum_{\sigma\in\mathcal{E}}\frac{1}{2}\beta^2|a^{w, n+1}_{K, \sigma}|^2|\sigma|d_{K, L}(w^{n+1}_{K} + w^{n+1}_{L})^2 + \frac{1}{2}\sum_{\sigma\in\mathcal{E}}\frac{|\sigma|}{d_{K,L}}(w^{n+1}_{K}-w^{n+1}_{L})^2 \nonumber \\
      &\leq 32\beta^2 C_J^2 \sum_{\sigma\in \mathcal{E}}|\sigma|d_{K, L} + \frac{1}{2}|w^{n+1}|^2_{1, \mathcal{T}} \leq C + \frac{1}{2}|w^{n+1}|^2_{1, \mathcal{T}}, \label{eq:drift_w_estimate} \\
      \sum_{K\in \mathcal{T}}\sum_{\sigma\in\mathcal{E}_K}\beta |\sigma|a^{v, n+1}_{K, \sigma}(v^{n+1}_{K} &+ v^{n+1}_{L})v^{n+1}_{K} = \sum_{\sigma\in\mathcal{E}}\beta |\sigma|\bigg\{a^{v, n+1}_{K, \sigma}(v^{n+1}_{K} + v^{n+1}_{L})v^{n+1}_{K} + a^{v, n+1}_{L, \sigma}(v^{n+1}_{K} + v^{n+1}_{L})v^{n+1}_{L} \bigg\} \nonumber \\
      &=\sum_{\sigma\in\mathcal{E}}\beta |\sigma|a^{v, n+1}_{K, \sigma}\bigg\{(v^{n+1}_{K} + v^{n+1}_{L})v^{n+1}_{K} -(v^{n+1}_{K} + v^{n+1}_{L})v^{n+1}_{L} \bigg\} \nonumber \\
      &=\sum_{\sigma\in\mathcal{E}}\beta |\sigma|a^{v, n+1}_{K, \sigma}\bigg(\big(v^{n+1}_K\big)^2 - \big(v^{n+1}_L\big)^2\bigg) \nonumber \\
      &\leq \sum_{\sigma\in\mathcal{E}}\beta |\sigma|\, |a^{v, n+1}_{K, \sigma}||v^{n+1}_{K} + v^{n+1}_{L}||v^{n+1}_{K}-v^{n+1}_{L}|\sqrt\frac{d_{K, L}}{d_{K,L}}  \nonumber \\
      &\leq \sum_{\sigma\in\mathcal{E}}\frac{1}{2}\beta^2|a^{v, n+1}_{K, \sigma}|^2|\sigma|d_{K, L}(v^{n+1}_{K} + v^{n+1}_{L})^2 + \frac{1}{2}\sum_{\sigma\in\mathcal{E}}\frac{|\sigma|}{d_{K,L}}(v^{n+1}_{K}-v^{n+1}_{L})^2 \nonumber \\
      &\leq 32\beta^2 C_J^2 \sum_{\sigma\in \mathcal{E}}|\sigma|d_{K, L} + \frac{1}{2}|v^{n+1}|^2_{1, \mathcal{T}} \leq C + \frac{1}{2}|v^{n+1}|^2_{1, \mathcal{T}}. \label{eq:drift_v_estimate}
\end{align}
To arrive at \eqref{eq:drift_w_estimate} and \eqref{eq:drift_v_estimate}, we have used that $|a^{\gamma, n+1}_{K, \sigma}|\leq C_J$ for $\gamma=\{w, v\}$ together with Young's inequality.
We also used that for a shape-regular quasi-uniform family of meshes in three
dimensions, one has $|\sigma|\leq Ch^2$ and $d_{K,L}\leq Ch$, while
$\#\sigma\leq Ch^{-3}$. Where $\#\sigma$ means the total number of faces in the domain. Consequently, it holds
\begin{equation*}
\sum_{\sigma\in\mathcal E}|\sigma|d_{K,L}
\leq Ch^2\cdot Ch\cdot Ch^{-3}
\leq C,
\end{equation*}
where $C$ is independent of $h$.\newline
Substituting the upper bounds obtained in \eqref{eq:dwdt_estimate}--\eqref{eq:drift_v_estimate} into the original equations for $w_K^{n}$ and $v^n_K$, we get

\begin{align}
    \frac{1}{2\tau}\Bigg(\sum_{K \in \mathcal{T}}\bigg\{|K|(w^{n+1}_K - w^{n}_K)^2 &+ |K|\big((w^{n+1}_K)^2 - (w^n_K)^2\big)\bigg\}\Bigg) +  |w^{n+1}|^2_{1, \mathcal{T}} + \sum_{\sigma\in\mathcal{E}} \beta |\sigma|\kappa_\sigma^{n+1} (w^{n+1}_K - w^{n+1}_{L})^2 \nonumber \\
    &- C - \frac{1}{2}|w^{n+1}|^2_{1, \mathcal{T}} \leq 0, \nonumber  \\
    \frac{1}{2\tau}\Bigg(\sum_{K \in \mathcal{T}}\bigg\{|K|(v^{n+1}_K - v^{n}_K)^2 &+ |K|\big((v^{n+1}_K)^2 - (v^n_K)^2\big)\bigg\}\Bigg) +  |v^{n+1}|^2_{1, \mathcal{T}} + \sum_{\sigma\in\mathcal{E}} \beta |\sigma|\kappa_\sigma^{n+1} (v^{n+1}_K - v^{n+1}_{L})^2 \nonumber \\
    &- C - \frac{1}{2}|v^{n+1}|^2_{1, \mathcal{T}} \leq 0. \nonumber 
\end{align}
Using the definition of our discrete norms \eqref{eq:H1norm_def}-\eqref{eq:H-1norm_def}, multiplying by $2\tau$, and then rearranging the resulting inequalities, we get 

\begin{align*}
    &\|w^{n+1}\|^2_{0, \mathcal{T}}-\|w^n\|^2_{0, \mathcal{T}} + \|w^{n+1}-w^{n}\|^2_{0, \mathcal{T}} + \tau |w^{n+1}|^2_{1, \mathcal{T}} + 2\tau\sum_{\sigma\in\mathcal{E}} \beta |\sigma|\kappa_\sigma^{n+1} (w^{n+1}_K - w^{n+1}_{L})^2 \leq  2\tau\,C, \\
    &\|v^{n+1}\|^2_{0, \mathcal{T}} - \|v^n\|^2_{0, \mathcal{T}} + \|v^{n+1}-v^{n}\|^2_{0, \mathcal{T}} + \tau |v^{n+1}|^2_{1, \mathcal{T}} + 2\tau \sum_{\sigma\in\mathcal{E}} \beta |\sigma|\kappa_\sigma^{n+1} (v^{n+1}_K - v^{n+1}_{L})^2 \leq  2\tau\,C. 
\end{align*}
Summing over $n\in\{0, \dots, N-1$\}, we see that 

\begin{align*}
    \sum_{n=0}^{N-1}\bigg(\|w^{n+1}\|^2_{0, \mathcal{T}} - \|w^{n}\|^2_{0, \mathcal{T}}\bigg) &= \|w^{N}\|_{0, \mathcal{T}} - \|w^{0}\|_{0, \mathcal{T}}, \\
    \sum_{n=0}^{N-1}\bigg(\|v^{n+1}\|^2_{0, \mathcal{T}} - \|v^{n}\|^2_{0, \mathcal{T}}\bigg) &= \|v^{N}\|_{0, \mathcal{T}} - \|v^{0}\|_{0, \mathcal{T}}.
\end{align*}
Consequently, we obtain the inequalities

\begin{align}
    &\|w^{N}\|^2_{0, \mathcal{T}} + \sum_{n=0}^{N-1}\|w^{n+1}-w^{n}\|^2_{0, \mathcal{T}} + \sum_{n=0}^{N-1}\tau |w^{n+1}|^2_{1, \mathcal{T}}  \leq \|w^0\|^2_{0, \mathcal{T}} + 2CT \leq C, \label{eq:normbounds_w_2} \\
    &\|v^{N}\|^2_{0, \mathcal{T}} + \sum_{n=0}^{N-1}\|v^{n+1}-v^{n}\|^2_{0, \mathcal{T}} + \sum_{n=0}^{N-1}\tau |v^{n+1}|^2_{1, \mathcal{T}}  \leq \|v^0\|^2_{0, \mathcal{T}} + 2CT \leq C. \label{eq:normbounds_v_2} 
\end{align}
In \eqref{eq:normbounds_w_2} and \eqref{eq:normbounds_v_2} we have removed the positive term stemming from the artificial diffusion. Note that $N-1$ is an arbitrary natural number and $C$ is independent of $h$ and $\tau$. Therefore, we arrive at the following estimates

\begin{align}
    &\sup_{0\leq n\leq N} \|w^n\|^2_{0, \mathcal{T}} + \sum_{n=0}^{N-1}\tau |w^{n+1}|^2_{1, \mathcal{T}}\leq C, \label{eq:Boch_H1_norm_w}\\
    &\sup_{0\leq n\leq N} \|v^n\|^2_{0, \mathcal{T}} + \sum_{n=0}^{N-1}\tau |v^{n+1}|^2_{1, \mathcal{T}}\leq C. \label{eq:Boch_H1_norm_v}
\end{align}
The inequalities in \eqref{eq:Boch_H1_norm_w} and \eqref{eq:Boch_H1_norm_v} leads to the first two inequalities stated in Lemma \ref{lem:EnergyBound}. 

\noindent\paragraph{Proof of \eqref{eq:Bochner_H-1_dt}} 
We introduce the notation $D_tw^{n+1}_K := \frac{w^{n+1}_K - w^n_K}{\tau}$ and $D_tv^{n+1}_K := \frac{v^{n+1}_K - v^n_K}{\tau}$ and write the discrete $H^{-1}$ norm as 

\begin{equation*}
    \|f\|_{-1, \mathcal{T}} = \sup_{\substack{
    \eta \in \mathcal{H}_{\mathcal T}\\
    \|\eta\|_{1,\mathcal T}\leq 1
}}\bigg\{\bigg|\sum_{K\in\mathcal{T}} |K| D_tf^{n+1}_K \eta_K\bigg|\bigg\}.
\end{equation*}

We fix $\eta \in \mathcal{H_T}$ such that $\|\eta\|_{1, \mathcal{T}}\leq 1$. Multiply \eqref{Eq:MainModelDisc_wv_w} and \eqref{Eq:MainModelDisc_wv_v} by $\eta$ to get

\begin{align}\label{eq:testingModelequations_w}
    \sum_{K\in\mathcal{T}}|K|D_tw^{n+1}_K\eta_K &= -\sum_{K\in\mathcal{T}}\sum_{\sigma\in\mathcal{E}_K}F^{w, n+1}_{K,\sigma}\eta_K = -\sum_{\sigma\in\mathcal{E}}F^{w, n+1}_{K,\sigma}\big(\eta_K -\eta_L \big), \\
    \sum_{K\in\mathcal{T}}|K|D_tv^{n+1}_K\eta_K &= -\sum_{K\in\mathcal{T}}\sum_{\sigma\in\mathcal{E}_K}F^{v, n+1}_{K,\sigma}\eta_K = -\sum_{\sigma\in\mathcal{E}}F^{v, n+1}_{K,\sigma}\big(\eta_K -\eta_L \big).\label{eq:testingModelequations_v}
\end{align}
We will treat each of the terms pertaining to the flux individually, starting with the diffusion. We have 
\begin{align}
    \bigg|\sum_{\sigma\in\mathcal{E}}\frac{|\sigma|}{d_{K,L}}(w^{n+1}_K-w^n_K)(\eta_K-\eta_L)\bigg| &\leq \Bigg(\frac{|\sigma|}{d_{K,L}}(w^{n+1}_K-w^n_K)^2\Bigg)^{1/2}\Bigg(\frac{|\sigma|}{d_{K,L}}(\eta_K-\eta_K)^2\Bigg)^{1/2} \nonumber \\
    &=|w^{n+1}|_{1, \mathcal{T}}\|\eta\|_{1, \mathcal{T}}, \label{eq:dualnorm_diffusion_w} \\
    \bigg|\sum_{\sigma\in\mathcal{E}}\frac{|\sigma|}{d_{K,L}}(v^{n+1}_K-v^n_K)(\eta_K-\eta_L)\bigg| &\leq \Bigg(\frac{|\sigma|}{d_{K,L}}(v^{n+1}_K-v^n_K)^2\Bigg)^{1/2}\Bigg(\frac{|\sigma|}{d_{K,L}}(\eta_K-\eta_K)^2\Bigg)^{1/2} \nonumber \\
    &=|v^{n+1}|_{1, \mathcal{T}}\|\eta\|_{1, \mathcal{T}}, \label{eq:dualnorm_diffusion_v}
\end{align}
and also
\begin{align}
     \bigg|\sum_{\sigma\in\mathcal{E}}\beta |\sigma|\kappa_{\sigma}^{n+1}(w^{n+1}_K-w^n_K)(\eta_K-\eta_L)\bigg| &\leq \beta \kappa_\sigma^{n+1}\Bigg(\frac{|\sigma|}{d_{K,L}}(w^{n+1}_K-w^n_K)^2\Bigg)^{1/2}\Bigg(|\sigma|d_{K, L}(\eta_K-\eta_K)^2\Bigg)^{1/2} \nonumber \\
    &\leq\beta \kappa_\sigma^{n+1}d_{K, L}|w^{n+1}|_{1, \mathcal{T}}|\eta|_{1, \mathcal{T}} \leq C_1|w^{n+1}|_{1, \mathcal{T}}\|\eta\|_{1, \mathcal{T}},\label{eq:dualnorm_numdiffusion_w} \\
    \bigg|\sum_{\sigma\in\mathcal{E}}\beta |\sigma|\kappa_{\sigma}^{n+1}(v^{n+1}_K-v^n_K)(\eta_K-\eta_L)\bigg| &\leq \beta \kappa_\sigma^{n+1}\Bigg(\frac{|\sigma|}{d_{K,L}}(v^{n+1}_K-v^n_K)^2\Bigg)^{1/2}\Bigg(|\sigma|d_{K, L}(\eta_K-\eta_K)^2\Bigg)^{1/2} \nonumber \\
    &\leq\beta \kappa_\sigma^{n+1}d_{K, L}|v^{n+1}|_{1, \mathcal{T}}|\eta|_{1, \mathcal{T}} \leq C_1|v^{n+1}|_{1, \mathcal{T}}\|\eta\|_{1, \mathcal{T}}.\label{eq:dualnorm_numdiffusion_v} 
\end{align}

Treating the drift terms results in 
\begin{align}
    \bigg|\sum_{\sigma \in \mathcal{E}}\beta |\sigma|a^{w, n+1}_{K, \sigma}(w^{n+1}_K + w^{n+1}_L)(\eta_K-\eta_L)\bigg| &\leq \sum_{\sigma\in \mathcal{E}}\frac{|\sigma|^{1/2}|\sigma|^{1/2} d_{K, L}^{1/2}}{d_{K, L}^{1/2}}|a^{w, n+1}_{K, \sigma}||w_K^{n+1}+w_L^{n+1}||\eta_K-\eta_L| \nonumber \\
    &\leq\beta |a^{w, n+1}_{K, \sigma}|\Bigg(\sum_{\sigma\in\mathcal{E}}|\sigma|d_{K, L}|w_K^{n+1} + w_L^{n+1}|^2\Bigg)^{1/2}\Bigg(\sum_{\sigma\in\mathcal{E}}\frac{|\sigma|}{d_{K,L}}|\eta_K-\eta_L|^2\Bigg)^{1/2} \nonumber \\
    &\leq 8\beta C_J \bigg(\sum_{\sigma\in\mathcal{E}}|\sigma|d_{K, L}\bigg)^{1/2}|\eta|_{1, \mathcal{T}} \leq C_2\|\eta\|_{1, \mathcal{T}}, \label{eq:dualnorm_drift_w} \\
    \bigg|\sum_{\sigma \in \mathcal{E}}\beta |\sigma|a^{v, n+1}_{K, \sigma}(v^{n+1}_K + v^{n+1}_L)(\eta_K-\eta_L)\bigg| &\leq \sum_{\sigma\in \mathcal{E}}\frac{|\sigma|^{1/2}|\sigma|^{1/2} d_{K, L}^{1/2}}{d_{K, L}^{1/2}}|a^{v, n+1}_{K, \sigma}||v_K^{n+1}+v_L^{n+1}||\eta_K-\eta_L| \nonumber \\
    &\leq\beta |a^{v, n+1}_{K, \sigma}|\Bigg(\sum_{\sigma\in\mathcal{E}}|\sigma|d_{K, L}|v_K^{n+1} + v_L^{n+1}|^2\Bigg)^{1/2}\Bigg(\sum_{\sigma\in\mathcal{E}}\frac{|\sigma|}{d_{K,L}}|\eta_K-\eta_L|^2\Bigg)^{1/2} \nonumber \\
    &\leq 8\beta C_J \bigg(\sum_{\sigma\in\mathcal{E}}|\sigma|d_{K, L}\bigg)^{1/2}|\eta|_{1, \mathcal{T}} \leq C_2\|\eta\|_{1, \mathcal{T}}. \label{eq:dualnorm_drift_v}
\end{align}

Gathering the estimates in
\eqref{eq:testingModelequations_w}--\eqref{eq:dualnorm_drift_v}, we obtain
\begin{align*}
    \left|\sum_{K\in\mathcal{T}} |K|D_tw^{n+1}_K\eta_K\right|
    &\leq \left((1+C_1)|w^{n+1}|_{1,\mathcal{T}}+C_2\right)
    \|\eta\|_{1,\mathcal{T}}, \\
    \left|\sum_{K\in\mathcal{T}} |K|D_tv^{n+1}_K\eta_K\right|
    &\leq \left((1+C_1)|v^{n+1}|_{1,\mathcal{T}}+C_2\right)
    \|\eta\|_{1,\mathcal{T}}. 
\end{align*}
Using $(1+a)^2\leq 2(1+a^2)$, we can find a constant $C>0$,
independent of $h$ and $\tau$, such that

\begin{align*}
    \left|\sum_{K\in\mathcal{T}} |K|D_tw^{n+1}_K\eta_K\right|^2
    &\leq C\left(1+|w^{n+1}|_{1,\mathcal{T}}^2\right)
    \|\eta\|_{1,\mathcal{T}}^2, \\
    \left|\sum_{K\in\mathcal{T}} |K|D_tv^{n+1}_K\eta_K\right|^2
    &\leq C\left(1+|v^{n+1}|_{1,\mathcal{T}}^2\right)
    \|\eta\|_{1,\mathcal{T}}^2. 
\end{align*}
Taking the supremum over all $\eta\in\mathcal{H}_{\mathcal{T}}$
such that $\|\eta\|_{1,\mathcal{T}}=1$, we obtain

\begin{align*}
    \|D_tw^{n+1}\|_{-1,\mathcal{T}}^2
    &\leq C\left(1+|w^{n+1}|_{1,\mathcal{T}}^2\right),\\
    \|D_tv^{n+1}\|_{-1,\mathcal{T}}^2
    &\leq C\left(1+|v^{n+1}|_{1,\mathcal{T}}^2\right).
\end{align*}
Multiplying the last two inequalities by $\tau$ and summing over $n=\{0,\ldots,N-1\}$ yields

\begin{align*}
    \sum_{n=0}^{N-1}\tau\,
    \|D_tw^{n+1}\|_{-1,\mathcal{T}}^2
    &\leq CT
    +C\sum_{n=0}^{N-1}\tau\,|w^{n+1}|_{1,\mathcal{T}}^2
    \leq C,\\
    \sum_{n=0}^{N-1}\tau\,
    \|D_tv^{n+1}\|_{-1,\mathcal{T}}^2
    &\leq CT
    +C\sum_{n=0}^{N-1}\tau\,|v^{n+1}|_{1,\mathcal{T}}^2
    \leq C,
\end{align*}
where we have used the $H^1$-estimates
\eqref{eq:Boch_H1_norm_w} and \eqref{eq:Boch_H1_norm_v}.
This proves the third inequality in Lemma ~\ref{lem:EnergyBound},
and hence, it completes the proof of the energy estimates for solutions of
Problem $(P_{h,\tau})$.
\end{proof}

\section{Well-posedness of approximate solutions to problem $(P_{h, \tau})$}\label{sec:well-posed}

Lemma \ref{lem:positivity} and Lemma \ref{lem:boundedness} are statements \emph{about solutions} of the scheme. This section shows that there exists at least one solution and later we prove that the solution is unique. The existence proof relies on a finite-dimensional fixed-point argument. We will fix the coefficients and solve the resulting linear system explicitly and then apply Brouwer's fixed-point theorem to prove that a solution exists. The set onto which the fixed-point is sought is a compact box determined by the initial mass of in the system and relies on the conservative flux, related to the calculations in Lemma \ref{lem:aprioribound}\footnote{Lemma \ref{lem:aprioribound} applies to the scheme with fixed parameters since it is still conserving mass.} and by the lower bound shown in Lemma \ref{lem:positivity}\footnote{Lemma \ref{lem:positivity} applies to the scheme with fixed parameters since the scheme retains the sign structure of the flux.}. The upper bound presented in Lemma \ref{lem:boundedness} is recovered \emph{a posteriori} at the fixed-point. The structure of our proof is similar to what is done in many other papers on finite volume schemes such as \cite{GAO202010} and \cite{FUHRMANN2001}. The argument to show that the system has a unique solution is by contradiction, we suppose that two solutions exists and show for sufficiently small time steps there difference is zero.

\subsection{Existence of solutions to problem $(P_{h, \tau})$}\label{sec:existence}
We prove the existence of at least one solution to problem $(P_{h, \tau})$, which is the main result of this section. The proof argues by induction and relies on Brouwer's Fixed-Point Theorem, restated below. For a proof of Brouwer's Fixed-Point Theorem, we direct the reader, for instance, to Section 3.2, Theorem 3.2 in \cite{Deimling1985}. 

\begin{theorem}[Brouwer's Fixed-Point Theorem] \label{thm:Brouwer}
    Let $X$ be a nonempty, compact and convex subset of a finite-dimensional normed vector space and let $f: X \to X$ be a continuous mapping. Then $f$ has at least one fixed-point.
\end{theorem}
Fix $h>0$ and define the set $$\mathcal{B}:=\{(\hat{w}, \hat{v})\in\mathbb{R}^{N_\mathcal{T}}\times\mathbb{R}^{N_\mathcal{T}}: 0\leq \hat{w}_K, \hat{v}_K \leq M_h \text{ for all }K\in\mathcal{T}\},$$ where $M_h$ is given in \eqref{eq:Mh}.
\begin{theorem}
    Assume $(H_1)$ and $(H_2)$. Problem $(P_{h,\tau})$ has at least one solution in the sense of Definition \ref{def:finite-volume-solution} lying in $\mathcal{B}$.
\end{theorem}
\begin{proof}
This proof proceeds by induction. Brouwer's Fixed-Point Theorem will be used to show that given a solution with the desired properties at the $n$ level (in the induction step), the solution map $\Phi$ has at least one fixed-point and that fixed-point satisfy problem $(P_{h, \tau})$ at the $n+1$ level.
\newline At $n=0$, we have $0\le |m_K^0|\le \phi_K^0\le 1$ and $w^0_K, v_K^0 \le M$ for any $K\in \mathcal{T}$, then it is clear that $w_K^0, v_K^0\ge 0$ and $w_K^0 + v^0_K=2\phi_K^0\le 2$. Since the flux is conservative, then we have 

\begin{align*}
    &\sum_{L\in \mathcal{T}}|L|w^0_K\le M|\Omega| \quad\text{ and }\quad \sum_{L\in \mathcal{T}}|L|v^0_K\le M|\Omega|.
\end{align*}
Fix $\hat{U}=(\hat{w}, \hat{v})$ $\in\mathcal{B}$. This pair is used to compute fixed coefficients using the formulas as anticipated earlier, viz.

\begin{align}
    \hat{m}_\sigma &:= \frac{\hat{w}_K-\hat{v}_K + \hat{w}_L - \hat{v}_L}{4},\label{eq:fixedCoeffs_1} \\
    \hat{V}_K &:= \sum_{L\in\mathcal{T}}|L|\nabla J(x_K-x_L)\frac{\hat{w}_L-\hat{v}_L}{2}, \\
    \hat{V}_\sigma &:= \frac{\hat{V}_K+\hat{V}_L}{2}, \\
    \hat{a}^w_{K, \sigma} &:=(1-\hat{m}_\sigma)\hat{V}_\sigma \cdot n_K,\\
    \hat{a}^v_{K, \sigma} &:=-(1+\hat{m}_\sigma)\hat{V}_\sigma \cdot n_K, \\
    \hat{\kappa}_\sigma &:= (1 + |\hat{m}_\sigma|)|\hat{V}_\sigma \cdot n_K|.\label{eq:fixedCoeffs_6}
\end{align}
Using these coefficients, we construct the flux dependent on $\hat{U}$ in terms of $(\hat{w}_K, \hat{v}_K)$:

\begin{align*}
    \hat{F}^{w}_{K, \sigma} &:=\frac{|\sigma|}{d_{K,L}}\big(\hat{w}_K-\hat{w}_L\big) + \beta |\sigma|\bigg(\big(\hat{a}^w_{K, \sigma} + \hat{\kappa}_\sigma\big)\hat{w}_K + \big(\hat{a}^w_{K, \sigma}-\hat{\kappa}_\sigma\big)\hat{w}_L\bigg) \\
    \hat{F}^{v}_{K, \sigma} &:=\frac{|\sigma|}{d_{K,L}}\big(\hat{v}_K-\hat{v}_L\big) + \beta |\sigma|\bigg(\big(\hat{a}^v_{K, \sigma} + \hat{\kappa}_\sigma\big)\hat{v}_K + \big(\hat{a}^v_{K, \sigma}-\hat{\kappa}_\sigma\big)\hat{v}_L\bigg) 
\end{align*}
Since $\hat{V}_\sigma$, $\hat{m}_\sigma$ are single valued on each face and $n_K=-n_L$, we see that $\hat{a}^w_{K, \sigma} = -\hat{a}^w_{L, \sigma}$ and $\hat{a}^v_{K, \sigma} = -\hat{a}^v_{L, \sigma}$. Consequently, this shows that $\hat{F}^w_{K, \sigma}=-\hat{F}^w_{L, \sigma}$ and 
$\hat{F}^v_{K, \sigma}=-\hat{F}^v_{L, \sigma}$. Given the $n$-level information $(w^n_K, v^n_K)$, we solve for $(\mathcal{W}_K, \mathcal{V}_K)$ for all $K\in \mathcal{T}$ the following system.

\begin{align}
      &\frac{|K|}{\tau} \mathcal{W}_K + \sum_{\sigma\in \mathcal{E}_K} \Bigg[\frac{|\sigma|}{d_{K,L}}(\mathcal{W}_K-\mathcal{W}_L) + \beta|\sigma|\bigg(\big(\hat{a}^w_{K, \sigma} + \hat{\kappa}_\sigma\big)\mathcal{W}_K + \big(\hat{a}^w_{K, \sigma} - \hat{\kappa}_\sigma\big)\mathcal{W}_L\bigg)\Bigg]=\frac{|K|}{\tau} w_K^n,
    \label{Eq:FrozenDisc_wv_w}
    \\[1ex]
    &\frac{|K|}{\tau} \mathcal{V}_K + \sum_{\sigma\in \mathcal{E}_K} \Bigg[\frac{|\sigma|}{d_{K,L}}( \mathcal{V}_K- \mathcal{V}_L) + \beta|\sigma|\bigg(\big(\hat{a}^v_{K, \sigma} + \hat{\kappa}_\sigma\big) \mathcal{V}_K + \big(\hat{a}^v_{K, \sigma} - \hat{\kappa}_\sigma\big) \mathcal{V}_L\bigg)\Bigg]=\frac{|K|}{\tau} v^n_K,
    \label{Eq:FrozenDisc_wv_v}
\end{align}
where the right-hand side only depends on the known information $(w^n_K, v_K^n)$. We rewrite \eqref{Eq:FrozenDisc_wv_w} and \eqref{Eq:FrozenDisc_wv_v} in matrix form to get
\begin{align}
    &A^w_{KK} \mathcal{W}_K + \sum_{L\in\mathcal{N}(K)}A^w_{KL}\mathcal{W}_L = \frac{|K|}{\tau}w^n_K,\label{eq:MatrixFrozen_w} \\
    &A^v_{KK} \mathcal{V}_K + \sum_{L\in\mathcal{N}(K)}A^v_{KL}\mathcal{V}_L = \frac{|K|}{\tau}v^n_K,\label{eq:MatrixFrozen_v} 
\end{align}
In \eqref{eq:MatrixFrozen_w}-\eqref{eq:MatrixFrozen_v}, the set $\mathcal{N}(K)$ denotes the neighbors of cell $K$ while the coefficients $A^w_{KK}$, $A^w_{KL}$, $A^v_{KK}$, and $A^v_{KL}$ are defined by:
\begin{align*}
    A^w_{KK}&:=\bigg(\frac{|K|}{\tau} + \sum_{\sigma\in\mathcal{E}_K}\bigg[\frac{|\sigma|}{d_{K,L}} + \beta|\sigma|\big(\hat{a}^w_{K, \sigma} + \hat{\kappa}_\sigma\big)\bigg]\bigg),&&A^w_{KL}:=-\frac{|\sigma|}{d_{K,L}} + \beta|\sigma|\big(\hat{a}^w_{K,\sigma} - \hat{\kappa}_\sigma\big),\\
    A^v_{KK}&:=\bigg(\frac{|K|}{\tau} + \sum_{\sigma\in\mathcal{E}_K}\bigg[\frac{|\sigma|}{d_{K,L}} + \beta|\sigma|\big(\hat{a}^v_{K, \sigma} + \hat{\kappa}_\sigma\big)\bigg]\bigg),&&A^v_{KL}:=-\frac{|\sigma|}{d_{K,L}} + \beta|\sigma|\big(\hat{a}^v_{K,\sigma} - \hat{\kappa}_\sigma\big).
\end{align*}
The above components define the matrices $A^w:=(A^w_{KL})_{K,L\in \mathcal{T}}$ and $A^v:=(A^v_{KL})_{K, L\in\mathcal{T}}$. We can write the solution map
$$\Phi:\mathcal B\to \mathcal B$$
via
\begin{equation*}
    \Phi(\hat{U}) := \bigg(\big((A^w(\hat{U}))^{-1}b^w\big), \big((A^v(\hat{U}))^{-1}b^v\big)\bigg) = \big(\mathcal{W}, \mathcal{V}\big),   \quad \text{with  }\,\mathcal{W}=(\mathcal{W_K})_{K\in\mathcal{T}}, \quad  \mathcal{V}= (\mathcal{V}_K)_{K\in\mathcal{T}}, 
\end{equation*}
where

\begin{align*}
    b^w&:= \bigg(\frac{|K|}{\tau}\hat{w}^n_K\bigg)_{K\in\mathcal{T}}, &&b^v := \bigg(\frac{|K|}{\tau}\hat{v}^n_K\bigg)_{K\in\mathcal{T}}.   
\end{align*}
To show that the system given by \eqref{eq:MatrixFrozen_w} and \eqref{eq:MatrixFrozen_v} has at least one solution, we show that the two matrices $A^w$ and $A^v$ are invertible and depend continuously on $(\hat{w}, \hat{v})$. We do so by demonstrating that these matrices are strictly diagonally dominant, see Theorem 5.6.15 and the associated Corollary 5.6.17 in \cite{HornJohnson2013} as well as the following calculations:

\begin{align*}
    A_{KK}^w &- \sum_{L\in\mathcal{N}(K)}|A_{LK}^w| \nonumber \\
    &=\frac{|K|}{\tau} + \sum_{\sigma\in\mathcal{E}_K}\bigg[\frac{|\sigma|}{d_{K,L}} + \beta|\sigma|\big(\hat{a}^w_{K, \sigma} + \hat{\kappa}_\sigma\big)\bigg] - \sum_{L\in\mathcal{N}(K)}\bigg|-\frac{|\sigma|}{d_{K,L}} + \beta|\sigma|\big(\hat{a}^w_{L,\sigma} - \hat{\kappa}_\sigma\big)\bigg| \nonumber \\
    &=\frac{|K|}{\tau} + \sum_{\sigma\in\mathcal{E}_K}\bigg[\frac{|\sigma|}{d_{K,L}} + \beta|\sigma|\big(\hat{a}^w_{K, \sigma} + \hat{\kappa}_\sigma\big)\bigg] - \sum_{\sigma\in\mathcal{E}_K}\bigg|-\frac{|\sigma|}{d_{K,L}} + \beta|\sigma|\big(-\hat{a}^w_{K,\sigma} - \hat{\kappa}_\sigma\big)\bigg| \nonumber \\
    &=\frac{|K|}{\tau} + \sum_{\sigma\in\mathcal{E}_K}\bigg[\frac{|\sigma|}{d_{K,L}} + \beta|\sigma|\big(\hat{a}^w_{K, \sigma} + \hat{\kappa}_\sigma\big)\bigg] - \sum_{\sigma\in\mathcal{E}_K}\frac{|\sigma|}{d_{K,L}} + \beta|\sigma|\big(\hat{a}^w_{K,\sigma} + \hat{\kappa}_\sigma\big)\nonumber \\
    &=\frac{|K|}{\tau} = \frac{h^3}{\tau} > 0. \label{eq:diagdom}
\end{align*}
The same expression holds for $A^v$. Therefore, both matrices $A^w$ and $A^v$ are strictly diagonally dominant independently of the fixed coefficients and thus invertible for every pair $(\hat{w}, \hat{v})\in \mathcal{B}$. The matrices $A^w$ and $A^v$ are continuous with respect to $\hat{U}$ as we can see by noticing that each element of $A^w$ and $A^v$ are affine transformations of $\hat{U}$, see \eqref{eq:fixedCoeffs_1}--\eqref{eq:fixedCoeffs_6}.
We still need to show that the wanted fixed point lies in $\mathcal{B}$. We preceed by induction, we assume at the $n$-level in the induction procedure that $(w^n_K, v_K^n)\in \mathcal{B}$, i.e. 
\begin{equation*}\label{eq:wn_vn_bounded}
\begin{aligned}
    &\sum_{L\in \mathcal{T}}|L|w^n_K\le M|\Omega|, &&\sum_{L\in \mathcal{T}}|L|v^n_K\le M|\Omega|,\\
    &w^n_K,v^n_K\ge 0, &&w^n_K, v_K^n\leq M_h.
\end{aligned}
\end{equation*}
For the candidate solutions $\mathcal{W}_K, \mathcal{V}_K$ obtained by solving \eqref{eq:MatrixFrozen_w} and \eqref{eq:MatrixFrozen_v}, we know by application of Lemma \ref{lem:positivity} and the conservative nature of the flux that

\begin{equation*}\label{eq:hat_boundedness}
\begin{aligned}
    &\sum_{K\in\mathcal{T}}|K|\mathcal{W}_K = \sum_{K\in \mathcal{T}}|K|w^n_K\leq M|\Omega|, && \sum_{K\in\mathcal{T}}|K|\mathcal{V}_K = \sum_{K\in \mathcal{T}}|K|v^n_K\leq M|\Omega|,\\
    &\mathcal{W}_K, \mathcal{V}_K\ge 0, &&\mathcal{W}_K, \mathcal{V}_K \le M_h.
\end{aligned}
\end{equation*}
 This ensures $(\mathcal{W}, \mathcal{V})\in \mathcal{B}$. Now, the hypothesis of Brouwer's Fixed-Point Theorem (restated in in Theorem \ref{thm:Brouwer}) are fulfilled and thus we know that there exists at least one $U^*\in \mathcal{B}$ such that $U^* = \Phi(U^*)$. At the fixed-point $U^*$, problem $(P_{h, \tau})$ is satisfied and thus we call $(w^{n+1}_K, v^{n+1}_K)=U^*$ the $n+1$ level solution. Since problem $(P_{h, \tau})$ is solved at the $n+1$ level, then we know that Lemma \ref{lem:boundedness} applies at the $n+1$ level and thus, we have $2\phi^{n+1}_K = w^{n+1}_K + v^{n+1}_K\le 2$.
\end{proof}

\subsection{Uniqueness of solutions to problem $(P_{h, \tau})$}\label{sec:uniqueness}
In this section we show that there is at most one solution $(w_K, v_K)$ to the system \eqref{Eq:MainModelDisc_wv}--\eqref{eq:fluxparams} in the sense of Definitions \ref{def:finite-volume-solution}. We begin by forming the equations for the difference of the suposed solutions and present some preliminaries before we state Theorem \ref{thm:uniqueness_induction} which is the main result of this section.
Let $(w^n_K, v^n_K)_{K\in\mathcal{T}}$ and $(\tilde{w}^n_K, \tilde{v}^n_K)_{K\in\mathcal{T}}$ be two distinct pairs of solutions to \eqref{Eq:MainModelDisc_wv}--\eqref{eq:fluxparams} with their accompanying parameters. Taking the difference between these sets of solutions, we obtain

\begin{align*}
    \frac{|K|}{\tau} (w_K^{n+1} - w_K^n)& - \frac{|K|}{\tau} (\tilde{w}_K^{n+1} - \tilde{w}_K^n) + \sum_{\sigma \in \mathcal{E}_K} \bigg(F^{w, n+1}_{K, \sigma} - F^{\tilde{w}, n+1}_{K, \sigma}\bigg) = 0, \\
    \frac{|K|}{\tau} (v_K^{n+1} - v_K^n)& - \frac{|K|}{\tau} (\tilde{v}_K^{n+1} - \tilde{v}_K^n) + \sum_{\sigma \in \mathcal{E}_K} \bigg(F^{v, n+1}_{K, \sigma} - F^{\tilde{v}, n+1}_{K, \sigma}\bigg) = 0.
\end{align*}
Let $e^w_K := w_K^{n+1} - \tilde{w}_K^{n+1}$, $e^v_K := v_K^{n+1} - \tilde{v}_K^{n+1}$ and we set $w_K^n = \tilde{w}_K^n$, $v_K^n = \tilde{v}_K^n$ since the solutions are assumed to be unique at $n$ in the following induction proof. This yields
\begin{subequations}\label{eq:diff_eq_wv}
\begin{align}
    \frac{|K|}{\tau}e^w_K &+ \sum_{\sigma \in \mathcal{E}_K} \bigg(F^{w}_{K, \sigma} - F^{\tilde{w}}_{K, \sigma}\bigg) = 0 , \\
    \frac{|K|}{\tau}e^v_K &+ \sum_{\sigma \in \mathcal{E}_K} \bigg(F^{v}_{K, \sigma} - F^{\tilde{v}}_{K, \sigma}\bigg) = 0 ,
\end{align}
\end{subequations}
where we have dropped the time-step index $n$ since every term is at the same time-step. Expanding \eqref{eq:diff_eq_wv}, we get 

\begin{align}
    \frac{|K|}{\tau } e^w_K &+ \sum_{\sigma \in \mathcal{E}_K}\frac{|\sigma|}{d_{K,L}}(e^w_K - e^w_L)\nonumber \\
    &+  \sum_{\sigma \in\mathcal{E}_K}\beta |\sigma|\bigg[(a^w_{K, \sigma} + \kappa_\sigma)w_K + (a^w_{K, \sigma} - \kappa_\sigma)w_L - (a^{\tilde{w}}_{K, \sigma} + \tilde{\kappa}_\sigma ) \tilde{w}_K - (a^{\tilde{w}}_{K, \sigma} - \tilde{\kappa }_\sigma )\tilde{w}_L\bigg] = 0, \label{eq:diff_eq_w}\\
    \frac{|K|}{\tau } e^v_K &+ \sum_{\sigma \in \mathcal{E}_K}\frac{|\sigma|}{d_{K,L}}(e^v_K - e^v_L)\nonumber \\
    &+  \sum_{\sigma \in\mathcal{E}_K}\beta |\sigma|\bigg[(a^v_{K, \sigma} + \kappa_\sigma)v_K + (a^v_{K, \sigma} - \kappa_\sigma)v_L - (a^{\tilde{v}}_{K, \sigma} + \tilde{\kappa}_\sigma ) \tilde{v}_K - (a^{\tilde{v}}_{K, \sigma} - \tilde{\kappa }_\sigma )\tilde{v}_L\bigg] = 0.\label{eq:diff_eq_v}
\end{align}
We add and subtract $(a^w_{K,\sigma} + \kappa_\sigma)\tilde{w}_K$ and $(a^w_{K, \sigma } - \kappa_\sigma)\tilde{w}_L$ in \eqref{eq:diff_eq_w}, while we add and subtract $(a^v_{K,\sigma} + \kappa_\sigma)\tilde{v}_K$ and $(a^v_{K, \sigma } - \kappa_\sigma)\tilde{v}_L$ in \eqref{eq:diff_eq_v}, to obtain 

\begin{align}
    &\frac{|K|}{\tau } e^w_K + \sum_{\sigma \in \mathcal{E}_K}\frac{|\sigma|}{d_{K,L}}(e^w_K - e^w_L)+  \sum_{\sigma \in\mathcal{E}_K}\beta |\sigma|\bigg[(a^w_{K, \sigma} + \kappa_\sigma)e^w_K + (a^w_{K, \sigma} - \kappa_\sigma)e^w_L\bigg]  = \rho^w_K, \label{eq:diffeq_w_rho}\\
    &\frac{|K|}{\tau } e^v_K + \sum_{\sigma \in \mathcal{E}_K}\frac{|\sigma|}{d_{K,L}}(e^v_K - e^v_L)+  \sum_{\sigma \in\mathcal{E}_K}\beta |\sigma|\bigg[(a^v_{K, \sigma} + \kappa_\sigma)e^v_K + (a^v_{K, \sigma} - \kappa_\sigma)e^v_L\bigg]  = \rho^v_K, \label{eq:diffeq_v_rho}
\end{align}
where 
\begin{align}
    &\rho_K^w := \sum_{\sigma\in \mathcal{E}_K} \beta |\sigma| \bigg[\big((a^{w}_{K, \sigma}-a^{\tilde{w}}_{K, \sigma}) + (\kappa_\sigma - \tilde{\kappa}_\sigma)\big)\tilde{w}_K + \big((a^{w}_{K, \sigma}-a^{\tilde{w}}_{K, \sigma}) + (\tilde{\kappa}_\sigma - \kappa_\sigma)\big)\tilde{w}_L\bigg], \label{eq:rho_w} \\
    &\rho_K^v := \sum_{\sigma\in \mathcal{E}_K} \beta |\sigma| \bigg[\big((a^{v}_{K, \sigma}-a^{\tilde{v}}_{K, \sigma}) + (\kappa_\sigma - \tilde{\kappa}_\sigma)\big)\tilde{v}_K + \big((a^{v}_{K, \sigma}-a^{\tilde{v}}_{K, \sigma}) + (\tilde{\kappa}_\sigma - \kappa_\sigma)\big)\tilde{v}_L\bigg]. \label{eq:rho_v}
\end{align}
For the purpose of clarity, we introduce the constants: $$C_J := \max_{K\in \mathcal{T}}\bigg\{\sum_{L \in \mathcal{T}}|L| |\nabla J (x_K-x_L)|\bigg\}$$ and $$\|e\|_\infty:= \max_{K\in\mathcal{T}}\max\bigg\{|e^w_K|, |e^v_K|\bigg\}.$$
We now show some key bounds that are direct consequences of Lemma \ref{lem:boundedness} and Lemma \ref{lem:boundedness}. They will turn out to be useful in estimating $\rho_K^w$ and $\rho_K^v$, viz.

\begin{align}
    |m_K-\tilde{m}_K| &= \bigg|\frac{w_K-v_K}{2} - \frac{\tilde{w}_K-\tilde{v}_K}{2}\bigg|=\frac{1}{2}|e^w_K - e^v_K|\leq \|e\|_\infty, \label{eq:mdiff_bound} \\
    \nonumber \\
    |V_K-\tilde{V}_K| &\leq \sum_{\sigma \in \mathcal{E}_L}|L||\nabla J(x_K-x_L)| |m_L-\tilde{m}_L| \leq C_J\|e\|_\infty, \label{eq:Vdiff_bound}\\
    \nonumber\\
    |a^{w}_{K, \sigma}-a^{\tilde{w}}_{K, \sigma}| &\leq \big|(1-m_\sigma)V_\sigma\cdot n_K - (1-\tilde{m}_\sigma) \tilde{V}_\sigma \cdot n_K\big| = \big|(1-m_\sigma)(V_\sigma - \tilde{V}_\sigma)\cdot n_K + (\tilde{m}_\sigma - m_\sigma)\tilde{V}_\sigma \cdot n_K \nonumber \\
    &\leq |1-m_\sigma| \, |(V_\sigma-\tilde{V}_\sigma )\cdot n_K| + |\tilde{m}_\sigma - m_\sigma|\,|\tilde{V}_\sigma \cdot n_K| \leq 3C_J\|e\|_\infty ,\label{eq:adiff_bound}\\
    \nonumber\\
    |\kappa_\sigma-\tilde{\kappa}_\sigma| & \leq \big|(1+|m_\sigma|)|V_\sigma\cdot n_K| - (1+|\tilde{m}_\sigma|) |\tilde{V}_\sigma \cdot n_K|\big|  \nonumber\\
    &=\big|(1+|m_\sigma|)(|V_\sigma\cdot n_K| - |\tilde{V}_\sigma\cdot n_K|) + (|\tilde{m}_\sigma| - |m_\sigma|)|\tilde{V}_\sigma \cdot n_K| \nonumber\\
    &\leq (1+|m_\sigma|) \, |(V_\sigma-\tilde{V}_\sigma )\cdot n_K| + |\tilde{m}_\sigma - m_\sigma|\,|\tilde{V}_\sigma \cdot n_K| \leq 3C_J\|e\|_\infty. \label{eq:kdiff_bound}
\end{align}
\newline
Now we are ready to formulate the uniqueness result. 
\begin{theorem}\label{thm:uniqueness_induction}
    Assume $(H_1)$-$(H_3)$. Then for every $n\in\mathbb{N}$, the problem $(P_{h, \tau})$ admits solutions in the sense of Definition \ref{def:finite-volume-solution} which is unique for all $K\in \mathcal{T}$.
\end{theorem}
\begin{proof}
The proof proceeds by induction. Clearly, uniqueness holds for $n=0$. The following discussion is showing that if uniqueness holds for a fixed arbitrary $n$, then it also holds for $n+1\in \mathbb{N}$. \newline
We apply the bounds in \eqref{eq:mdiff_bound}, \eqref{eq:Vdiff_bound}, \eqref{eq:adiff_bound}, and \eqref{eq:kdiff_bound} to $\rho_K^w$ and $\rho_K^v$ in \eqref{eq:rho_w} and \eqref{eq:rho_v}, respectively, to get

\begin{align*}
    |\rho_K^w| &\leq \sum_{\sigma\in \mathcal{E}_K} \beta |\sigma| \bigg[\big(|a^{w}_{K, \sigma}-a^{\tilde{w}}_{K, \sigma}| + |\kappa_\sigma - \tilde{\kappa}_\sigma|\big)\tilde{w}_K + \big(|a^{w}_{K, \sigma}-a^{\tilde{w}}_{K, \sigma}| + |\tilde{\kappa}_\sigma - \kappa_\sigma| \big)\tilde{w}_L\bigg] \nonumber \\
    &\leq \sum_{\sigma\in \mathcal{E}_K} \beta |\sigma|  \bigg[\big(3C_J \|e\|_\infty + 3C_J \|e\|_\infty\big)2 +  \big(3C_J \|e\|_\infty + 3C_J \|e\|_\infty\big)2\bigg] \nonumber \\
    & \leq 24\beta C_J \|e\|_\infty\sum_{\sigma \in \mathcal{E}_K}|\sigma| \leq 144\beta C_J h^2\|e\|_\infty, \label{eq:unique_upper_bound_w} \\
    |\rho_K^v| & \leq \sum_{\sigma\in \mathcal{E}_K} \beta |\sigma| \bigg[\big(|a^{v}_{K, \sigma}-a^{\tilde{v}}_{K, \sigma}| + |\kappa_\sigma - \tilde{\kappa}_\sigma|\big)\tilde{v}_K + \big(|a^{v}_{K, \sigma}-a^{\tilde{v}}_{K, \sigma}| + |\tilde{\kappa}_\sigma - \kappa_\sigma| \big)\tilde{v}_L\bigg] \nonumber \\
    &\leq 24 \beta C_J \|e\|_\infty\sum_{\sigma \in \mathcal{E}_K}|\sigma| \leq 144\beta C_J h^2\|e\|_\infty,
\end{align*}
where the last inequality follows from the six identical sides of each cell.
Turning our attention to the left-hand side of \eqref{eq:diffeq_w_rho}, let $K^*$ be a certain $K\in\mathcal{T}$ where $e^w_K$ attains its maximum. We assume, without loss of generality, that $w_{K^*}\geq \tilde{w}_{K^*}$. For the chosen $K^*$, it holds that $e^w_{K^*}\geq e^w_L$. We then have

\begin{align*}
    \frac{|K^*|}{\tau}\|e\|_\infty &+ \sum_{\sigma \in\mathcal{E}_{K^*}}\bigg[\frac{|\sigma|}{d_{K,L}}(\|e\|_\infty-e^w_L)  + \beta |\sigma | \bigg(\big(a^w_{K, \sigma} + \kappa_\sigma\big)\|e\|_\infty + \big(a^w_{K, \sigma} - \kappa_\sigma\big)e^w_L  \bigg)\bigg] \\
    &\geq \frac{|K^*|}{\tau}\|e\|_\infty + \beta\sum_{\sigma \in \mathcal{E}_K}|\sigma|\bigg(a^w_{K, \sigma}\big(\|e\|_\infty + e^w_L\big) + \kappa_\sigma\big(\|e\|_\infty - e^w_L\big)\bigg) \\
    &\geq\frac{|K^*|}{\tau}\|e\|_\infty +\beta \sum_{\sigma \in \mathcal{E}_K}|\sigma|a^{w}_{K, \sigma}\big(\|e\|_\infty + e^w_L\big).
\end{align*}
To move on, we use that if $a^w_{K, \sigma}\geq0$ then $a^w_{K, \sigma}(\|e\|_\infty + e^w_L)\geq 0$, together with, if $a^w_{K, \sigma}<0$,  then $a^w_{K, \sigma}(\|e\|_\infty + e^w_L)\geq -2|a^w_{K, \sigma}|\|e\|_\infty$. We get

\begin{align*}
     \frac{|K^*|}{\tau}\|e\|_\infty &+ \sum_{\sigma \in\mathcal{E}_K}\frac{|\sigma|}{d_{K,L}}(\|e\|_\infty-e^w_L)  + \beta |\sigma | \bigg(\big(a^w_{K, \sigma} + \kappa_\sigma\big)\|e\|_\infty + \big(a^w_{K, \sigma} - \kappa_\sigma\big)e^w_L  \bigg) \nonumber  \\
    &\geq \frac{|K^*|}{\tau }\|e\|_\infty -\sum_{\sigma\in\mathcal{E}_K}2\beta|\sigma||a^w_{K, \sigma}|\,\|e\|_\infty \geq \frac{|K^*|}{\tau }\|e\|_\infty - 24\beta C_J h^2 \|e\|_\infty. 
\end{align*}
Following the same reasoning applied to \eqref{eq:diffeq_v_rho}, let $\tilde{K}^*$ be a certain $K\in \mathcal{T}$ where $e^v_K$ attains its maximum assuming $v_{\tilde{K}^*}\geq \tilde{v}_{\tilde{K}^*}$. For the chosen $\tilde{K}^*$, it holds that $e^v_{\tilde{K}^*}\geq e^v_L$. We get 
\begin{align*}
    \frac{|K^*|}{\tau}\|e\|_\infty &+ \sum_{\sigma \in\mathcal{E}_K}\bigg[\frac{|\sigma|}{d_{K,L}}(\|e\|_\infty-e^v_L)  + \beta |\sigma | \bigg(\big(a^v_{K, \sigma} + \kappa_\sigma\big)\|e\|_\infty + \big(a^w_{K, \sigma} - \kappa_\sigma\big)e^v_L  \bigg)\bigg] \nonumber\\
    &\geq \frac{|K^*|}{\tau}\|e\|_\infty + \beta\sum_{\sigma \in \mathcal{E}_K}|\sigma|\bigg(a^v_{K, \sigma}\big(\|e\|_\infty + e^v_L\big) + \kappa_\sigma\big(\|e\|_\infty - e^v_L\big)\bigg) \nonumber\\
    &\geq\frac{|K^*|}{\tau}\|e\|_\infty +\beta \sum_{\sigma \in \mathcal{E}_K}|\sigma|a^{v}_{K, \sigma}\big(\|e\|_\infty + e^v_L\big) \geq \frac{|K^*|}{\tau }\|e\|_\infty - 24\beta C_J h^2 \|e\|_\infty. 
\end{align*}
Collecting the upper and lower bounds for both $e^w_K$ and $e^v_K$, we finally obtain the following result for each equation

\begin{equation*}
    (\frac{h^3}{\tau} - 172\beta C_J h^2 ) \|e\|_\infty \leq 0.
\end{equation*}
Pick $\tau_0:=\dfrac{h}{172\beta C_J}$ and let  $\tau < \tau_0$. This allows us conclude that $\|e\|_\infty = 0$ through

\begin{equation*}
    \Big(\tau_0 - \tau\Big)\|e\|_\infty \leq 0, \label{eq:dt_restiction}
\end{equation*}
which finally proves uniqueness at the $n+1$ level, thus completing the induction argument.
\end{proof}

\begin{remark}
    The solution to Problem $(P_{h, \tau})$ is stable with respect to parameters and initial data. The proof of this statement follows the same style of argument as in the proof of uniqueness of solutions done for Theorem \ref{thm:uniqueness_induction}. For a proof of stability with respect to parameters in the semi-discrete setting for the same model, see Theorem 3.1 in our previous work \cite{NicklasJ2026}.
\end{remark}
\section{Compactness result and identification of the limit equations} \label{sec:compactness}
Let $(\mathcal{T}_i, \tau_i)_{i\in \mathbb{N}}$ be a family of discretizations, where $\mathcal{T}_i$ is an admissible mesh on $\Omega$ with $h_i = \underset{K\in\mathcal{T}_i}{\max}\text{diam}(K)$ such that $h_i\to 0$ while $i\to \infty$. Let $D_t^{\tau_i} w^{n+1}: = \frac{w^{n+1}-w^n}{\tau_i}$ and $N_i:=\frac{T}{\tau_i}$. 
For convenience, in Lemma \ref{discrete-lions-aubin-verbatim} we restate Lemma 9 from \cite{Cances2019} and then apply it to problem $(P_{h, \tau})$. Lemma \ref{discrete-lions-aubin-verbatim} is the key ingredients for passing to the limit in our approximate solutions.
\begin{lemma}[Discrete Lions-Aubin]\label{discrete-lions-aubin-verbatim}
    Let $(u_i)\subset \mathcal{H}_{\mathcal{T}_i, \tau_i}$ be a sequence of piecewise constant in time functions with values in $\mathcal{H}_{\mathcal{T}_i}$ satisfying 
    \begin{align}
        &\sum_{n=0}^{N_i}\bigg( \tau \|u^{n}_i\|^2_{1, \mathcal{T}_i} + \tau\|D_t^{\tau_i}u_i^{n}\|^2_{-1, \mathcal{T}_m}\bigg) \leq C, \label{eq:u_energy}
    \end{align}
    where $C> 0$ is independent of the size of the mesh and the time step size. Then there exists a subsequence, which is not relabeled, such that, as $i\to \infty$, it holds
    \begin{align*}
        u_i &\to u \qquad \text{  strongly in }L^2(Q_T),\\
        \nabla u_i &\rightharpoonup \nabla u \qquad \text{weakly  in }L^2(Q_T).
    \end{align*}
\end{lemma}

\noindent Applying Lemma \ref{discrete-lions-aubin-verbatim} to our scheme, we see that \eqref{eq:u_energy} is satisfied for both $(w^{n}_K)$ and $(v^n_K)$. Thus we can find convergent subsequences as stated in Theorem \ref{thm:compactness}.
\begin{theorem} \label{thm:compactness}
Assume $(H_1)$-$(H_3)$. Let $(w_i), (v_i)\subset \mathcal{H}_{\mathcal{T}_i, \tau_i}$ be sequences of piecewise constant in time functions with values in $\mathcal{H}_{\mathcal{T}_i}$ solutions to problem $(P_{h, \tau})$ in the sense of Definition \ref{def:finite-volume-solution}, satisfying the \emph{a priori} bounds

\begin{align*}
    &\sum_{n=0}^{N_i-1}\bigg( \tau \|w^{n+1}_i\|^2_{1, \mathcal{T}_i} + \tau\|D_t^{\tau_i}w_i^{n+1}\|^2_{-1, \mathcal{T}_i}\bigg) \leq C, \\
    &\sum_{n=0}^{N_i-1}\bigg( \tau \|v^{n+1}_i\|^2_{1, \mathcal{T}_i} + \tau\|D_t^{\tau_i}v^{n+1}_i\|^2_{-1, \mathcal{T}_i} \bigg)\leq C,
\end{align*}
where $C> 0$ is independent of the size of the mesh and of the time step size. Then there exists a subsequences, which are not relabeled, such that, as $i\to \infty$, it holds

\begin{align*}
    w_i&\to w \qquad \text{strongly in } L^2(Q_T),\\
    v_i&\to v \,\qquad \text{strongly in } L^2(Q_T), \\
    \nabla w_i&\rightharpoonup \nabla w \,\quad \text{weakly in } L^2(Q_T),\\
    \nabla v_i&\rightharpoonup \nabla v \,\,\quad \text{weakly in } L^2(Q_T). 
\end{align*}
\end{theorem}
\begin{proof}
    The conclusion of this Theorem is a direct consequence of Lemma \ref{lem:EnergyBound} combined with Lemma \ref{discrete-lions-aubin-verbatim}.
\end{proof}

\begin{corollary} \label{cor:compactness}
    Recalling $m_i = \frac{w_i-v_i}{2}$ and $\phi_i = \frac{w_i + v_i}{2}$, we know that 
    \begin{align*}
    m_i&\to m \qquad \text{strongly in } L^2(Q_T),\\
    \phi_i&\to \phi \,\qquad \text{strongly in } L^2(Q_T), \\
    \nabla m_i&\rightharpoonup \nabla m \,\quad \text{weakly in } L^2(Q_T),\\
    \nabla \phi_i&\rightharpoonup \nabla \phi \,\,\quad \text{weakly in } L^2(Q_T). 
\end{align*}
as $i\to \infty$.
\end{corollary}

\section{Identifying the limit}\label{sec:limit}

\begin{theorem}
    Assume $(H_1)$-$(H_3)$. Then the limit points identified in Theorem \ref{thm:compactness} arrived at by solving problem $(P_{h, \tau})$ are weak solutions in the sense of Definition \ref{Def:MainModelWeak_wv}.
\end{theorem}
\begin{proof}
In this section, we pass to the limit $i\to \infty$ and show that the discrete scheme given in \eqref{Eq:MainModelDisc_wv}--\eqref{eq:fluxparams} converges to the weak solution of \eqref{Eq:MainModel_wv}. Let $w_{\mathcal{T}_i, \tau_i}, v_{\mathcal{T}_i, \tau_i}\in \mathcal{H}_{\mathcal{T}_i, \tau_i}$,  where $\mathcal{H}_{\mathcal{T}_i, \tau_i}$ is given in \eqref{eq:H_T_t} with $\mathcal{T}=\mathcal{T}_i$ and $\tau =\tau_i$ with $h=\sup\{\text{diam} (K):K\in \mathcal{T}_i\}$. 
Take $\varphi\in C^\infty_0(\Bar{\Omega}\times[0, T))$. We multiply the governing equations in \eqref{Eq:MainModelDisc_wv} by $\tau\varphi(t^n, x_K)=\tau \varphi^n_K$ and sum over all $K\in \mathcal{T}$ and all $n\in \{0, \dots, N-1\}$ to obtain

\begin{align*}
    T_1^w + T_2^w + T_3^w + T_4^w = 0, \\
    T_1^v + T_2^v + T_3^v + T_4^v = 0, 
\end{align*}
where 
\begin{align}
    &T_1^w :=\sum_{n=0}^{N-1}\sum_{K\in\mathcal{T}}|K|(w^{n+1}_K-w^{n}_K)\varphi^n_K, &&T_1^v :=\sum_{n=0}^{N-1}\sum_{K\in\mathcal{T}}|K|(v^{n+1}_K-v^{n}_K)\varphi^n_K, \nonumber \\
    &T_2^w := \sum_{n=0}^{N-1}\tau\sum_{K\in\mathcal{T}}\sum_{\sigma\in \mathcal{E}_K}\frac{|\sigma|}{d_{K,L}}(w^{n+1}_K-w^{n+1}_L)\varphi^n_K, &&T_2^v := \sum_{n=0}^{N-1}\tau\sum_{K\in\mathcal{T}}\sum_{\sigma\in \mathcal{E}_K}\frac{|\sigma|}{d_{K,L}}(v^{n+1}_K-v^{n+1}_L)\varphi^n_K, \nonumber \\
    &T_3^w :=\sum_{n=0}^{N-1}\tau\sum_{K\in\mathcal{T}}\sum_{\sigma\in \mathcal{E}_K}\beta |\sigma|\kappa_\sigma^{n+1}(w^{n+1}_K+w^{n+1}_L)\varphi^n_K, &&T_3^v :=\sum_{n=0}^{N-1}\tau\sum_{K\in\mathcal{T}}\sum_{\sigma\in \mathcal{E}_K}\beta |\sigma|\kappa_\sigma^{n+1}(v^{n+1}_K+v^{n+1}_L)\varphi^n_K, \nonumber \\
    &T_4^w :=\sum_{n=0}^{N-1}\tau\sum_{K\in\mathcal{T}}\sum_{\sigma\in \mathcal{E}_K}\beta |\sigma|a^{w,n+1}_{K, \sigma}(w^{n+1}_K+w^{n+1}_L)\varphi^n_K, &&T_4^v :=\sum_{n=0}^{N-1}\tau\sum_{K\in\mathcal{T}}\sum_{\sigma\in \mathcal{E}_K}\beta |\sigma|a^{v,n+1}_{K, \sigma}(v^{n+1}_K+v^{n+1}_L)\varphi^n_K. \nonumber 
\end{align}
It is important to remember that $\varphi$ is a smooth function evaluated at discrete points, thus allowing us to utilize Taylor's expansion formula to rewrite 
\begin{align}
    \varphi^n_L =\varphi(t^n, x_L) &= \varphi(t^n, x_K) + (x_L-x_K)\nabla\varphi(t^n, x_K) + \mathcal{O}(d^2_{K, L})\nonumber \\
    &=\varphi(t^n, x_K) + d_{K, L}\nabla\varphi(t^n, x_K)\cdot n_K + \mathcal{O}(d^2_{K, L})\label{eq:Taylor}
\end{align}
Note that if $h\to 0$, then $d_{K,L}\to 0$. \\
We will now treat $T_j^{\gamma}$ for $j\in\{1, 2, 3, 4\}$ and $\gamma\in \{w, v\}$, separately. 
The discrete integration by-parts-in-time applied to $T_1^{w}$ and $T_1^v$ results in 

\begin{align*}
    T_1^w &= \sum_{n=1}^{N-1} \sum_{K\in \mathcal{T}}|K|w^n_K(\varphi^{n-1}_K-\varphi^n_K) + \sum_{K\in\mathcal{T}}|K|w^N_K\varphi^{N-1}_K -\sum_{K\in\mathcal{T}}|K|w^0_K\varphi^0_K \nonumber \\
    &=-\sum_{n=0}^{N-1}\sum_{K\in \mathcal{T}}\int_{t^{n-1}}^{t^n}\int_K w^n_K\partial_t\varphi(t, x_K)dxdt, \\
    T_1^v &= \sum_{n=1}^{N-1} \sum_{K\in \mathcal{T}}|K|v^n_K(\varphi^{n-1}_K-\varphi^n_K) + \sum_{K\in\mathcal{T}}|K|v^N_K\varphi^{N-1}_K -\sum_{K\in\mathcal{T}}|K|v^0_K\varphi^0_K\nonumber \\
    &= -\sum_{n=0}^{N-1}\sum_{K\in\mathcal{T}}\int_{t^{n-1}}^{t^n}\int_K v^n_K\partial_t\varphi(t, x_K)dxdt,
\end{align*}
where we have used $\varphi(T-\tau, x_K)\to 0$ as $\tau\to 0$. This allows us pass to the limit to obtain
\begin{align}
    T_1^w\to &-\int_0^T\int_\Omega w\partial_t\varphi dxdt \text{ as } i\to \infty \nonumber, \\
    T_1^v\to &-\int_0^T\int_\Omega v\partial_t\varphi dxdt \text{ as } i\to \infty. \nonumber 
\end{align}
Moving on to treating $T_2^w$ and $T_2^v$, we rewrite them such that we sum over each face once, i.e.

\begin{align}
    T_2^w &= \sum_{n=0}^{N-1}\tau\sum_{\sigma\in \mathcal{E}}|\sigma|(w^{n+1}_K-w^{n+1}_L)\frac{\varphi^n_K-\varphi^n_L}{d_{K,L}}, \nonumber \\
    T_2^v &= \sum_{n=0}^{N-1}\tau\sum_{\sigma\in \mathcal{E}}|\sigma|(v^{n+1}_K-v^{n+1}_L)\frac{\varphi^n_K-\varphi^n_L}{d_{K,L}}. \nonumber 
\end{align}
At this point we introduce two terms, $\tilde{T}^w_2$ and $\tilde{T}^v_2$, that converge to the desired limits, $\tilde{T}^w_2$ is derived as follows

\begin{align}
    \tilde{T}^w_2 :&=\sum_{n=0}^{N-1}\int_{n\tau}^{(n+1)\tau}\int_\Omega w_{\mathcal{T}, \tau}(x, t)\Delta \varphi(x, t^n) dx\,dt\, \nonumber \\
    &=\sum_{n=0}^{N-1}\tau \sum_{K\in\mathcal{T}}w^n_K\int_K \Delta \varphi(x,t^n) \,dx\,dt \nonumber \\
    &=\sum_{n=0}^{N-1}\tau \sum_{K\in\mathcal{T}}w^n_K\int_{\partial K}\nabla \varphi(x, t^n)\cdot n_K \,d\sigma \nonumber \\
    &=\sum_{n=0}^{N-1}\tau \sum_{K\in\mathcal{T}}\sum_{\sigma\in\mathcal{E}_K}w^n_K \int_\sigma\nabla \varphi(x, t^n)\cdot n_K \,d\sigma \nonumber \\
    &=\sum_{n=0}^{N-1}\tau \sum_{\sigma \in \mathcal{E}}(w^n_K-w^n_L)\int_\sigma \nabla \varphi(x, t^n)\cdot n_K \,d\sigma. \nonumber
\end{align}
Similarly, we have for $\tilde{T}^v_2$ the following expression

\begin{equation*}
    \tilde{T}^v_2 := \sum_{n=0}^{N-1}\tau \sum_{\sigma \in \mathcal{E}}(v^n_K-v^n_L)\int_\sigma \nabla \varphi(x, t^n)\cdot n_K \,d\sigma. 
\end{equation*}
The terms $T^w_2 + \tilde{T}^w_2$ and $T^w_2 + \tilde{T}^w_2$ can now be written in the following way 

\begin{align*}
    T^w_2 + \tilde{T}^w_2 &= \sum_{n=0}^{N-1}\tau \sum_{\sigma \in \mathcal{E}}|\sigma|(w^{n+1}_K-w^{n+1}_L)R^{n}_{K,L}, \\
    T^v_2 + \tilde{T}^v_2 &= \sum_{n=0}^{N-1}\tau \sum_{\sigma \in \mathcal{E}}|\sigma|(v^{n+1}_K-v^{n+1}_L)R^{n}_{K,L},
\end{align*}
where we get the remainder

\begin{align*}
    R^{n}_{K,L} &= \frac{1}{|\sigma|}\int_\sigma \nabla \varphi(x, t^n) \cdot n_K\,d\sigma - \frac{\varphi^n_L - \varphi^n_L}{d_{K, L}}.
\end{align*}
If we use \eqref{eq:Taylor} to expand $\varphi$ around a point on the surface $\sigma$, referred to as $y_\sigma$, to express both $\varphi^n_L$ and $\varphi^n_K$, then we get 

\begin{equation*}
    R^n_{K, L} =\frac{1}{|\sigma|}\int_\sigma \nabla \varphi(x, t^n)\cdot n_K\,d\sigma -\nabla \varphi(y_\sigma, t^n)\cdot n_K + \mathcal{O}(d_{K, L}).
\end{equation*}
We take the absolute value and then utilize the Mean-Value Theorem to obtain

\begin{align*}
    |R^n_{K, L}| \leq \frac{1}{|\sigma|}\int_\sigma |\nabla \varphi(x, t^n) - \nabla \varphi(y_\sigma, t^n)|\,d\sigma \leq \frac{1}{|\sigma|}\int_\sigma  C |\frac{d}{d\xi} \nabla \varphi(\xi, t^n)| \,h\,d\sigma \leq C\|\frac{d}{d\xi}\nabla \varphi\|_{L^\infty(\Omega)} h, 
\end{align*}
where $\xi$ is some point on the straight line segment joining $x$ and $y_\sigma$. This argument shows that

\begin{align*}
    T^w_2 &\to -\int_\Omega w(x,t) \Delta \varphi(x, t)\, dx \quad \text{ as } i \to \infty, \\
    T^v_2 &\to -\int_\Omega v(x,t) \Delta \varphi(x, t)\, dx \quad \text{ as } i \to \infty.
\end{align*}
Anticipating the passing to the limit for the artificial diffusion, notice that the artificial diffusive terms, $T^w_3$ and $T^v_3$, are of order $h^2$ whilst the original diffusion in $T^w_2$ and $T^v_2$ is on the order of $h$. Therefore, we should expect that the artificial diffusion will tend to zero as the discretization becomes finer. 

\begin{align}
    T_3^w &=\sum_{n=0}^{N-1}\tau \sum_{\sigma \in \mathcal{E}}\beta |\sigma|\kappa^{n+1}_\sigma (w^{n+1}_K-w^{n+1}_L)(\varphi^n_K-\varphi^n_L) \nonumber \\
    &=\sum_{n=0}^{N-1}\tau \sum_{\sigma \in \mathcal{E}}\beta |\sigma|\kappa^{n+1}_\sigma (w^{n+1}_K-w^{n+1}_L)(d_{K,L}\nabla \varphi^n_K\cdot n_K + \mathcal{O}(d^{2}_{K,L})) \nonumber \\
    &\leq 8\beta C_J|\Omega|\sum_{n=0}^{N-1}\tau \sum_{\sigma\in\mathcal{E}}d_{K, L}\Big(\int_\sigma |\nabla \varphi^n_K|d\sigma + \mathcal{O}(d^2_{K, L})\Big)\nonumber \\
    &\leq 8\beta C_J|\Omega|\sum_{n=0}^{N-1}\tau \sum_{\sigma\in\mathcal{E}}h\Big(\int_\sigma |\nabla \varphi^n_K|d\sigma + \mathcal{O}(d^2_{K, L})\Big)\label{eq:T3_w}.
\end{align}
A similar expression is obtained for $T^v_3$, i.e.

\begin{align}
    T^v_3 &\leq 8\beta C_J|\Omega|\sum_{n=0}^{N-1}\tau \sum_{\sigma\in\mathcal{E}}h\Big(\int_\sigma |\nabla \varphi^n_K|d\sigma + \mathcal{O}(d^2_{K, L})\Big),  \label{eq:T3_v}
\end{align}
To obtain \eqref{eq:T3_w} and \eqref{eq:T3_v}, we have used the upper bound $\kappa_\sigma^{n+1}\leq 2C_J|\Omega|$. Since $|\nabla \varphi|$ is bounded from above, it forces $T^w_3$ and $T^v_3$ to vanish as $h\to 0$.\\
Finally, we rewrite $T_4^w$ and $T^v_4$ using the definition of $a^{w, n+1}_{K, \sigma}$ and $a^{v, n+1}_{K, \sigma}$, viz.

\begin{align}
    T_4^w &= \sum_{n=0}^{N-1}\tau\sum_{\sigma\in \mathcal{E}}\beta |\sigma|(1-m_\sigma^{n+1})(w_K^{n+1}+ w_L^{n+1})(\varphi_K^n - \varphi_L^n)V_\sigma^{n+1}\cdot n_K \nonumber \\
    &=-\frac{\beta}{2}\sum_{n=0}^{N-1}\tau\sum_{K\in \mathcal{T}}\sum_{\sigma\in \mathcal{E}_K} |\sigma|(1-m^{n+1}_\sigma)(w_K^{n+1}+ w_L^{n+1})(d_{K,L}\nabla \varphi_\sigma^{n}\cdot n_K + \mathcal{O}(d^2_{K,L}))V_\sigma^{n+1}\cdot n_K, \label{eq:T4_w_partial}  \\
    T_4^v &=\frac{\beta}{2}\sum_{n=0}^{N-1}\tau\sum_{K\in \mathcal{T}}\sum_{\sigma\in \mathcal{E}_K} |\sigma|(1+m^{n+1}_\sigma)(v_K^{n+1}+ v_L^{n+1})(d_{K,L}\nabla \varphi_\sigma^{n}\cdot n_K + \mathcal{O}(d^2_{K,L}))V_\sigma^{n+1}\cdot n_K. \label{eq:T4_v_partial}
\end{align}
To obtain \eqref{eq:T4_w_partial} and \eqref{eq:T4_v_partial}, we applied \eqref{eq:Taylor} to $\varphi^n_K$ and $\varphi^n_L$ and denoted $\nabla \varphi_\sigma^n =\frac{1}{2}\big(\nabla\varphi_K^n + \nabla \varphi_L^n\big)$.
In dealing with \eqref{eq:T4_w_partial} and \eqref{eq:T4_v_partial}, we propose two terms, $\Tilde{T}^w_4$ and $\tilde{T}^v_4$, respectively. We show that $|T^w_4-\Tilde{T}^w_4|\to 0$ and $|T^v_4-\tilde{T}^v_4|\to 0$ as $h\to 0$. we show that $\tilde{T}^w_4$ and $\tilde{T}^v_4$ will converge to the expected continuum limit.
The terms $\tilde{T}^w_4$ and $\tilde{T}^v_4$ are:

\begin{align}
    \tilde{T}^w_4 :&=-2\beta \sum_{n=0}^{N-1}\tau \sum_{K\in\mathcal{T}}|K|(1-m^{n+1}_K)w^{n+1}_KV^{n+1}_K\cdot \nabla \varphi_K^n \nonumber \\
    &= -\frac{\beta}{2}\sum_{n=0}^{N-1}\tau \sum_{K\in\mathcal{T}}\sum_{\sigma\in \mathcal{E}_K}2|\sigma|d_{K,L}(1-m_K^{n+1})w^{n+1}_K(V_K\cdot n_K)(\nabla \varphi^n_K\cdot n_K). \label{eq:T4_w*} \\
    \tilde{T}^v_4 :&= \frac{\beta}{2}\sum_{n=0}^{N-1}\tau \sum_{K\in\mathcal{T}}\sum_{\sigma\in \mathcal{E}_K}2|\sigma|d_{K,L}(1+m_K^{n+1})v^{n+1}_K(V_K\cdot n_K)(\nabla \varphi^n_K\cdot n_K). \label{eq:T4_v*}
\end{align}
To derive the structure of \eqref{eq:T4_w*} and \eqref{eq:T4_v*}, we utilized the fact that for any two vectors $a,b\in \mathbb{R}^n$ we have 

\begin{equation*}
    \sum_{\sigma \in\mathcal{E}_K}|\sigma|d_{K,L}(a\cdot n_K)(b\cdot n_K) =2 |K| (a\cdot b).
\end{equation*}
We proceed by defining $D^w_{K, \sigma}$ in the following manner

\begin{align}
    D^w_{K, \sigma} &= (1-m^{n+1}_\sigma)(w^{n+1}_K+w^{n+1}_L)(V_\sigma^{n+1}\cdot n_K)(\nabla \varphi^{n}_\sigma\cdot n_K) - 2(1-m^{n+1}_K)w^{n+1}_K(V^{n+1}_K\cdot n_K)(\nabla \varphi^n_K\cdot n_K) \nonumber \\
    &=(1-m^{n+1}_\sigma)(w^{n+1}_K+w^{n+1}_L)(V_\sigma^{n+1}\cdot n_K)(\nabla (\varphi^{n}_\sigma-\nabla \varphi_K^n)\cdot n_K) \nonumber\\
    &+(1-m^{n+1}_\sigma)(w^{n+1}_K+w^{n+1}_L)(V_\sigma^{n+1}\cdot n_K)(\nabla \varphi^{n}_K\cdot n_K) \nonumber \\
    &-2(1-m^{n+1}_K)w^{n+1}_K(V^{n+1}_K\cdot n_K)(\nabla \varphi^n_K\cdot n_K) \nonumber \\
    &=(1-m^{n+1}_\sigma)(w^{n+1}_K+w^{n+1}_L)(V_\sigma^{n+1}\cdot n_K)(\nabla (\varphi^{n}_\sigma-\nabla \varphi_K^n)\cdot n_K) \nonumber \\
    &+(1-m^{n+1}_\sigma)(w^{n+1}_K+w^{n+1}_L)((V_\sigma^{n+1}-V_K^{n+1})\cdot n_K)(\nabla \varphi^{n}_K\cdot n_K) \nonumber \\
    &+(1-m^{n+1}_\sigma)(w^{n+1}_K+w^{n+1}_L)(V_K^{n+1}\cdot n_K)(\nabla \varphi^{n}_K\cdot n_K) \nonumber \\
    &-2(1-m^{n+1}_K)w^{n+1}_K(V^{n+1}_K\cdot n_K)(\nabla \varphi^n_K\cdot n_K) \nonumber \\
    &=(1-m^{n+1}_\sigma)(w^{n+1}_K+w^{n+1}_L)(V_\sigma^{n+1}\cdot n_K)(\nabla (\varphi^{n}_\sigma-\nabla \varphi_K^n)\cdot n_K) \nonumber\\
    &+(1-m^{n+1}_\sigma)(w^{n+1}_K+w^{n+1}_L)((V_\sigma^{n+1}-V_K^{n+1})\cdot n_K)(\nabla \varphi^{n}_K\cdot n_K)  \nonumber\\
    &+(1-m^{n+1}_\sigma)(w^{n+1}_L-w^{n+1}_K)(V_K\cdot n_K)(\nabla \varphi^n_K\cdot n_K) \nonumber\\
    &+ 2w^{n+1}_K(m^{n+1}_K-m^{n+1}_\sigma)(V_K\cdot n_K)(\nabla \varphi^n_K\cdot n_K) .\nonumber
\end{align}
A similar expression holds for $D^v_{K, \sigma}$, namely

\begin{align*}
     D^v_{K, \sigma} &=(1+m^{n+1}_\sigma)(v^{n+1}_K+v^{n+1}_L)(V_\sigma^{n+1}\cdot n_K)(\nabla (\varphi^{n}_\sigma-\nabla \varphi_K^n)\cdot n_K) \\
    &+(1+m^{n+1}_\sigma)(v^{n+1}_K+v^{n+1}_L)((V_\sigma^{n+1}-V_K^{n+1})\cdot n_K)(\nabla \varphi^{n}_K\cdot n_K) \\
    &+(1+m^{n+1}_\sigma)(v^{n+1}_L-v^{n+1}_K)(V_K\cdot n_K)(\nabla \varphi^n_K\cdot n_K)  \\
    &+ 2v^{n+1}_K(m^{n+1}_\sigma-m^{n+1}_K)(V_K\cdot n_K)(\nabla \varphi^n_K\cdot n_K) .
\end{align*}
For convenience in writing, we introduce the next level of notation. 

\begin{align*}
    E^{1, w}_{K, \sigma}:&=(1-m^{n+1}_\sigma)(w^{n+1}_K+w^{n+1}_L)(V_\sigma^{n+1}\cdot n_K)(\nabla (\varphi^{n}_\sigma-\nabla \varphi_K^n)\cdot n_K)  \\
    E^{2, w}_{K,\sigma}:&=(1-m^{n+1}_\sigma)(w^{n+1}_K+w^{n+1}_L)((V_\sigma^{n+1}-V_K^{n+1})\cdot n_K)(\nabla \varphi^{n}_K\cdot n_K)
    \\
    E^{3, w}_{K, \sigma}:&=(1-m^{n+1}_\sigma)(w^{n+1}_L-w^{n+1}_K)(V_K\cdot n_K)(\nabla \varphi^n_K\cdot n_K) 
    \\
    E^{4, w}_{K, \sigma}:&= 2w^{n+1}_K(m^{n+1}_K-m^{n+1}_\sigma)(V_K\cdot n_K)(\nabla \varphi^n_K\cdot n_K), 
\end{align*}
and 

\begin{align*}
    E^{1, v}_{K, \sigma}:&=(1+m^{n+1}_\sigma)(v^{n+1}_K+v^{n+1}_L)(V_\sigma^{n+1}\cdot n_K)(\nabla (\varphi^{n}_\sigma-\nabla \varphi_K^n)\cdot n_K)  \\
    E^{2, v}_{K,\sigma}:&=(1+m^{n+1}_\sigma)(v^{n+1}_K+v^{n+1}_L)((V_\sigma^{n+1}-V_K^{n+1})\cdot n_K)(\nabla \varphi^{n}_K\cdot n_K)
    \\
    E^{3, v}_{K, \sigma}:&=(1+m^{n+1}_\sigma)(v^{n+1}_L-v^{n+1}_K)(V_K\cdot n_K)(\nabla \varphi^n_K\cdot n_K) 
    \\
    E^{4, v}_{K, \sigma}:&= 2v^{n+1}_K(m^{n+1}_\sigma-m^{n+1}_K)(V_K\cdot n_K)(\nabla \varphi^n_K\cdot n_K). 
\end{align*}
Dropping the term that scales as $d^2_{K,L}$ in \eqref{eq:T4_w_partial}, we write

\begin{align*}
    T^w_4-\tilde{T}^w_4 &= -\frac{\beta}{2}\sum_{n=0}^{N-1}\tau\sum_{K\in\mathcal{T}}\sum_{\sigma\in\mathcal{E}_K}|\sigma|d_{K,L}D^w_{K, \sigma} \nonumber \\
    &= -\frac{\beta}{2}\sum_{n=0}^{N-1}\tau\sum_{K\in\mathcal{T}}\sum_{\sigma\in\mathcal{E}_K}|\sigma|d_{K,L}\big(E^{1,w}_{K, \sigma}+E^{2,w}_{K, \sigma}+E^{3,w}_{K, \sigma}+E^{4,w}_{K, \sigma}\big).
\end{align*}
Similarly, we drop the $\mathcal{O}(d^2_{K,L})$-term in \eqref{eq:T4_v_partial} to write
\begin{align*}
    T^v_4-\tilde{T}^v_4 &= \frac{\beta}{2}\sum_{n=0}^{N-1}\tau\sum_{K\in\mathcal{T}}\sum_{\sigma\in\mathcal{E}_K}|\sigma|d_{K,L}D^v_{K, \sigma} \nonumber \\
    &= \frac{\beta}{2}\sum_{n=0}^{N-1}\tau\sum_{K\in\mathcal{T}}\sum_{\sigma\in\mathcal{E}_K}|\sigma|d_{K,L}\big(E^{1,v}_{K, \sigma}+E^{2,v}_{K, \sigma}+E^{3,v}_{K, \sigma}+E^{4,v}_{K, \sigma}\big).
\end{align*}
We now deal separately with each term of the sum. We do this explicitly for $T^w_4-\tilde{T}^w_4$, and state the corresponding result for $T^v_4-\tilde{T}^v_4$.

\begin{align}
    -\frac{\beta}{2}\sum_{n=0}^{N-1}\tau& \sum_{K\in\mathcal{T}}\sum_{\sigma\in\mathcal{E}_K}|\sigma|d_{K,L}E^{1,w}_{K,\sigma} \nonumber \\
    &= -\frac{\beta}{2}\sum_{n=0}^{N-1}\tau \sum_{K\in\mathcal{T}}\sum_{\sigma\in\mathcal{E}_K}|\sigma|d_{K,L}(1-m^{n+1}_\sigma)(w^{n+1}_K+w^{n+1}_L)(V^{n+1}_\sigma\cdot n_K)((\nabla \varphi_\sigma^n-\nabla \varphi^n_K)\cdot n_K) \nonumber \\
    &\le \frac{\beta}{2}\sum_{n=0}^{N-1}\tau \sum_{K\in\mathcal{T}}\sum_{\sigma\in\mathcal{E}_K}|\sigma|d_{K,L}|1-m^{n+1}_\sigma|\,|w^{n+1}_K+w^{n+1}_L|\,|V^{n+1}_\sigma|\,|\nabla \varphi_\sigma^n-\nabla \varphi^n_K| \nonumber \\
    &\leq 2\beta |\Omega||\nabla J|_{L^\infty(\Omega)}\sum_{n=0}^{N-1}\tau \sum_{K\in\mathcal{T}}\sum_{\sigma\in\mathcal{E}_K}|\sigma|d_{K,L}|\nabla \varphi_L^n-\nabla \varphi^n_K| \nonumber \\
    &\le 2\beta |\Omega||\|\nabla J\|_{L^\infty(\Omega)}\|\sum_{n=0}^{N-1}\tau \sum_{K\in\mathcal{T}}\sum_{\sigma\in\mathcal{E}_K}|\sigma|d_{K,L}\big(|d_{K, L}\Delta \varphi^n_K| + \mathcal{O}(d^2_{K,L})\big)\nonumber \\
    &\le 2\beta |\Omega|\|\nabla J\|_{L^\infty(\Omega)}\|\Delta\varphi\|_{L^\infty(\Omega)} h \sum_{n=0}^{N-1}\tau \sum_{K\in\mathcal{T}}\sum_{\sigma\in\mathcal{E}_K}|\sigma|d_{K,L} + \mathcal{O}(d^2_{K,L})\nonumber \\
    &\le 6\beta |\Omega|^2\|\nabla J\|_{L^\infty(\Omega)}\|\Delta\varphi\|_{L^\infty(\Omega)} hT + \mathcal{O}(d^2_{K,L}). \nonumber 
\end{align}
Thus, we have
\begin{align}
    -\frac{\beta}{2}\sum_{n=0}^{N-1}\tau& \sum_{K\in\mathcal{T}}\sum_{\sigma\in\mathcal{E}_K}|\sigma|d_{K,L}E^{1,w}_{K,\sigma} \leq 6\beta |\Omega|^2\|\nabla J\|_{L^\infty(\Omega)}\|\Delta\varphi\|_{L^\infty(\Omega)} hT + \mathcal{O}(d^2_{K,L})
    \label{eq:E1_w_lim}, \\
    \frac{\beta}{2}\sum_{n=0}^{N-1}\tau& \sum_{K\in\mathcal{T}}\sum_{\sigma\in\mathcal{E}_K}|\sigma|d_{K,L}E^{1,v}_{K,\sigma} \le 6\beta |\Omega|^2\|\nabla J\|_{L^\infty(\Omega)}\|\Delta\varphi\|_{L^{\infty}(\Omega)}h T + \mathcal{O}(d^2_{K,L}) \label{eq:E1_v_lim}.
\end{align}
Looking at \eqref{eq:E1_w_lim} and \eqref{eq:E1_v_lim} we conclude that the $(E^{1,w}_{K, \sigma})$ and $(E^{1,v}_{K, \sigma})$-terms tends to zero as $h\to 0$. Continuing with $E^{2,\cdot}_{K, \sigma}$ we write

\begin{align}
    -\frac{\beta}{2}\sum_{n=0}^{N-1}\tau& \sum_{K\in\mathcal{T}}\sum_{\sigma\in\mathcal{E}_K}|\sigma|d_{K,L}E^{2,w}_{K,\sigma} \nonumber \\
    &=-\frac{\beta}{2}\sum_{n=0}^{N-1}\tau\sum_{K\in\mathcal{T}}\sum_{\sigma\in\mathcal{E}_K}|\sigma|d_{K,L}(1-m^{n+1}_\sigma)(w^{n+1}_K+w^{n+1}_L)((V^{n+1}_\sigma-V^{n+1}_K)\cdot n_K)(\nabla \varphi^n_K\cdot n_K) \nonumber \\
    &\le \frac{\beta}{2}\sum_{n=0}^{N-1}\tau \sum_{K\in\mathcal{T}}\sum_{\sigma\in\mathcal{E}_K}|\sigma|d_{K,L}|1-m^{n+1_\sigma}|\,|w^{n+1}_K+w^{n+1}_L|\,|V^{n+1}_\sigma-V^{n+1}_K|\,|\nabla \varphi^n_K| \nonumber \\
    &\le 2\beta \|\nabla \varphi\|_{L^\infty(\Omega)}\sum_{n=0}^{N-1}\tau\sum_{K\in\mathcal{T}}\sum_{\sigma\in\mathcal{E}_K}|\sigma|d_{K,L}|V^{n+1}_L-V^{n+1}_K|\nonumber \\
    &\le 2\beta|\Omega| \|\nabla \varphi\|_{L^\infty(\Omega)}\|\Delta J\|_{L^\infty(\Omega)}h\sum_{n=0}^{N-1}\tau \sum_{K\in\mathcal{T}}\sum_{\sigma\in\mathcal{E}_K}|\sigma|d_{K,L} + \mathcal{O}(d^2_{K,L})\nonumber. 
\end{align}
Thus, we conclude 

\begin{align}
    -\frac{\beta}{2}\sum_{n=0}^{N-1}\tau \sum_{K\in\mathcal{T}}\sum_{\sigma\in\mathcal{E}_K}|\sigma|d_{K,L}E^{2,w}_{K,\sigma}&\le 6\beta|\Omega|^2 \|\nabla \varphi\|_{L^\infty(\Omega)}\|\Delta J\|_{L^\infty(\Omega)}hT+\mathcal{O}(d^2_{K,L}), \label{eq:E2_w_lim} \\
    \frac{\beta}{2}\sum_{n=0}^{N-1}\tau \sum_{K\in\mathcal{T}}\sum_{\sigma\in\mathcal{E}_K}|\sigma|d_{K,L}E^{2,v}_{K,\sigma} &\le  6\beta|\Omega|^2 \|\nabla \varphi\|_{L^\infty(\Omega)}\|\nabla^2J\|_{L^\infty(\Omega)}hT + \mathcal{O}(d^2_{K,L}). \label{eq:E2_v_lim}
\end{align}
This proves that the left-hand sides in \eqref{eq:E2_w_lim} and \eqref{eq:E2_v_lim} tend to zero as $h\to 0$. In dealing with the last two terms, recall the definition of the discrete gradient norm in \eqref{eq:H1norm_def}.
We proceed as follows:
\begin{align}
    -\frac{\beta}{2}\sum_{n=0}^{N-1}\tau& \sum_{K\in\mathcal{T}}\sum_{\sigma\in\mathcal{E}_K}|\sigma|d_{K,L}E^{3,w}_{K,\sigma} \nonumber \\
    &= -\frac{\beta}{2}\sum_{n=0}^{N-1}\tau \sum_{K\in\mathcal{T}}\sum_{\sigma\in\mathcal{E}_K}|\sigma|d_{K,L}(1-m^{n+1}_\sigma)(w^{n+1}_L-w^{n+1}_K)(V^{n+1}_\sigma\cdot n_K)(\nabla \varphi^n_K\cdot n_K) \nonumber \\
    &\le \frac{\beta}{2}\sum_{n=0}^{N-1}\tau \sum_{K\in\mathcal{T}}\sum_{\sigma\in\mathcal{E}_K}|\sigma|d_{K,L}|1-m^{n+1}_\sigma|\,|w^{n+1}_L-w^{n+1}_K|\,|V^{n+1}_\sigma|\,|\nabla \varphi^n_K| \nonumber \\
    &\le \beta |\Omega|\|\nabla J\|_{L^\infty(\Omega)}\|\nabla \varphi\|_{L^\infty(\Omega)}\sum_{n=0}^{N-1}\tau \sum_{K\in\mathcal{T}}\sum_{\sigma\in\mathcal{E}_K}|\sigma|^{1/2}d^{3/2}_{K,L}\frac{|\sigma|^{1/2}}{d^{1/2}_{K,L}}|w^{n+1}_L-w^{n+1}_K| \nonumber \\
    & \le 2\beta |\Omega|\|\nabla J\|_{L^\infty(\Omega)}\|\nabla \varphi\|_{L^\infty(\Omega)}\sum_{n=0}^{N-1}\tau \Bigg(\sum_{\sigma\in\mathcal{E}} |\sigma| d^3_{K, L}\Bigg)^{1/2}\Bigg(\sum_{\sigma\in\mathcal{E}} \frac{|\sigma|}{d_{K, L}}|w^{n+1}_L-w^{n+1}_K|^2\Bigg)^{1/2} \nonumber \\
    &\le 2\sqrt{3}\beta |\Omega|^{3/2} \|\nabla J\|_{L^\infty(\Omega)}\|\nabla \varphi\|_{L^\infty(\Omega)} h \sum_{n=0}^{N-1}\tau |w^{n+1}|_{1, \mathcal{T}}, \nonumber 
\end{align}
where we have used the Cauchy-Schwarz inequality. Summarizing, we get 

\begin{align}
    -\frac{\beta}{2}\sum_{n=0}^{N-1}\tau \sum_{K\in\mathcal{T}}\sum_{\sigma\in\mathcal{E}_K}|\sigma|d_{K,L}E^{3,w}_{K,\sigma} &\leq 2\sqrt{3}\beta |\Omega|^{3/2} \|\nabla J\|_{L^\infty(\Omega)}\|\nabla \varphi\|_{L^\infty(\Omega)} h \sum_{n=0}^{N-1}\tau |w^{n+1}|_{1, \mathcal{T}}, \label{eq:E3_w_lim}\\
    \frac{\beta}{2}\sum_{n=0}^{N-1}\tau \sum_{K\in\mathcal{T}}\sum_{\sigma\in\mathcal{E}_K}|\sigma|d_{K,L}E^{3,v}_{K,\sigma} &\le 2\sqrt{3}\beta |\Omega|^{3/2} \|\nabla J\|_{L^\infty(\Omega)}\|\nabla \varphi\|_{L^\infty(\Omega)} h \sum_{n=0}^{N-1}\tau |v^{n+1}|_{1, \mathcal{T}}, \label{eq:E3_v_lim}
\end{align}
where we recall \eqref{eq:Boch_H1_norm_w} and \eqref{eq:Boch_H1_norm_v} to determine that \eqref{eq:E3_w_lim} and \eqref{eq:E3_v_lim} tends to zero as $h\to 0$. We turn our attention to the final term

\begin{align}
    -\frac{\beta}{2}\sum_{n=0}^{N-1}\tau& \sum_{K\in\mathcal{T}}\sum_{\sigma\in\mathcal{E}_K}|\sigma|d_{K,L}E^{4,w}_{K,\sigma} \nonumber \\
    &=-\beta\sum_{n=0}^{N-1}\tau \sum_{K\in\mathcal{T}}\sum_{\sigma\in\mathcal{E}_K}|\sigma|d_{K,L}|w_k|(m^{n+1}_k-m^{n+1}_\sigma)(V^{n+1}_K\cdot n_K)(\nabla \varphi^n_K\cdot n_K) \nonumber \\
    &\le \beta\sum_{n=0}^{N-1}\tau \sum_{K\in\mathcal{T}}\sum_{\sigma\in\mathcal{E}_K}|\sigma|d_{K,L}|w_k|\,|m^{n+1}_k-m^{n+1}_\sigma|\,|V^{n+1}_K|\,|\nabla \varphi^n_K| \nonumber \\
    &\le 2\beta|\Omega| \|\nabla J\|_{L^\infty(\Omega)}\|\nabla \varphi\|_{L^\infty(\Omega)}\sum_{n=0}^{N-1}\tau \sum_{K\in\mathcal{T}}\sum_{\sigma\in\mathcal{E}_K}|\sigma|^{1/2}d^{3/2}_{K,L}\frac{|\sigma|^{1/2}}{d^{1/2}_{K,L}}|m^{n+1}_K-m^{n+1}_L|\nonumber \\
    &\le 4\beta|\Omega| \|\nabla J\|_{L^\infty(\Omega)}\|\nabla \varphi\|_{L^\infty(\Omega)} \sum_{n=0}^{N-1}\tau\Bigg(\sum_{\sigma\in\mathcal{E}}|\sigma|d_{K,L}^3\Bigg)^{1/2}\Bigg(\sum_{\sigma \in \mathcal{E}}\frac{|\sigma |}{d_{K, L}}|m^{n+1}_K-m^{n+1}_L|\Bigg)^{1/2} \nonumber \\
    &\le 4\sqrt{3}|\Omega|^{3/2}\|\nabla J\|_{L^\infty(\Omega)}\|\nabla \varphi\|_{L^\infty(\Omega)} h\sum_{n=0}^{N-1}\tau |m^{n+1}|_{1, \mathcal{T}} \nonumber \\
    &\le 4\sqrt{3}|\Omega|^{3/2}\|\nabla J\|_{L^\infty(\Omega)}\|\nabla \varphi\|_{L^\infty(\Omega)} h\sum_{n=0}^{N-1}\tau \big(|w^{n+1}|_{1, \mathcal{T}} + |v^{n+1}|_{1, \mathcal{T}}\big). \nonumber 
\end{align}
Finally, this allows us to write 
\begin{align*}
    -\frac{\beta}{2}\sum_{n=0}^{N-1}\tau \sum_{K\in\mathcal{T}}\sum_{\sigma\in\mathcal{E}_K}|\sigma|d_{K,L}E^{4,w}_{K,\sigma}&
    \le 4\sqrt{3}|\Omega|^{3/2}\|\nabla J\|_{L^\infty(\Omega)}\|\nabla \varphi\|_{L^\infty(\Omega)} h\sum_{n=0}^{N-1}\tau \big(|w^{n+1}|_{1, \mathcal{T}} + |v^{n+1}|_{1, \mathcal{T}}\big), 
    \\
    \frac{\beta}{2}\sum_{n=0}^{N-1}\tau \sum_{K\in\mathcal{T}}\sum_{\sigma\in\mathcal{E}_K}|\sigma|d_{K,L}E^{4,v}_{K,\sigma} &\le4\sqrt{3}|\Omega|^{3/2}\|\nabla J\|_{L^\infty(\Omega)}\|\nabla \varphi\|_{L^\infty(\Omega)} h\sum_{n=0}^{N-1}\tau \big(|w^{n+1}|_{1, \mathcal{T}} + |v^{n+1}|_{1, \mathcal{T}}\big). 
\end{align*}
With the knowledge that $|T^w_4-\tilde{T}^w_4|\to 0$ and $|T^v_4-\tilde{T}^v_4|\to 0$ as $h\to 0$, we are allowed to use \eqref{eq:T4_w*} and \eqref{eq:T4_v*} to show to what limit our scheme converge.
Passing to the limit in \eqref{eq:T4_w*} and \eqref{eq:T4_v*} results in 

\begin{align*}
    &\tilde{T}^w_4 \to -2\beta \int_\Omega (1-m)w\nabla J\star m\nabla \varphi \, dx\quad \text{ as }i \to \infty, 
    &\tilde{T}^v_4 \to 2\beta \int_\Omega (1-m)w\nabla J\star m\nabla \varphi \, dx\quad \text{ as }i \to \infty. 
\end{align*}
This concludes our proof that the solutions of the proposed numerical scheme converge to the solutions of the weak formulation presented in Definition \ref{Def:MainModelWeak_wv}.
\end{proof}

\section*{Acknowledgments} NJ and AM were involved in Swedish Energy Agency’s project Solar Electricity Research Centre
(SOLVE) with grant number 52693-1. 
AM also acknowledges the Swedish Research Council for partial financial support under project  no. 2024-05606.



\bibliographystyle{plain}

\end{document}